\documentclass{article}
\usepackage{amssymb}
\usepackage{graphicx}
\usepackage{amsmath}

\input{tcilatex}
\begin{document}

\begin{center}
\textbf{Multipliers and} \textbf{embedding operators with application to
abstract differential equat\i ons }

{\textbf{Veli B. Shakhmurov}}

Okan University, Department of Mechanical engineering, Akfirat, Tuzla 34959
Istanbul, Turkey,

E-mail: veli.sahmurov@okan.edu.tr

\bigskip\ \ \ \ \ \ \ \ \ \ \ \ \ \ 
\end{center}

\QTP{Body Math}
\ \ \ \ \ \ \ \ \ \ \ \ \ \ \ \ \ \ \ \ \ \ \ \ \ \ \ \ \ \ \ \ \ \ \ \ \ \
\ \ \ \ \ \ \ \ \ \ \ \ \ \ 

\begin{center}
\textbf{ABSTRACT}
\end{center}

\begin{quote}
\ \ \ \ \ \ \ \ \ \ \ \ \ \ \ 
\end{quote}

\ \ In this paper,\ Mikhlin and Marcinkiewicz--Lizorkin type operator-valued
multiplier theorems in weighted Lebesgue-Bochner spaces are studied. Using
these results one derives embedding theorems in $E_{0}$-valued weighted
Sobolev-Lions type spaces $W_{p,\gamma }^{l}\left( \Omega ;E_{0},E\right) $,
where $E_{0}$, $E$ are two Banach spaces, $E_{0}$ is continuously and
densely embedded into $E.$ One proves that, there exists a smoothest
interpolation space $E_{\alpha },$ between $E_{0}$ and $E$, such that the
differential operator $D^{\alpha }$ acts as a bounded linear operator from $%
W_{p,\gamma }^{l}\left( \Omega ;E_{0},E\right) $ to $L_{p,\gamma }\left(
\Omega ;E_{\alpha }\right) $. By using these results the $L_{p,\gamma }-$%
separability properties of elliptic operators and regularity properties of
appropriate degenerate differential operators are studied. In particular, we
prove that the associated differential operator is positive and also is a
generator of an analytic semigroup. Moreover, the maximal $L_{p,\gamma }$%
-regularity properties of Cauchy problem for abstract parabolic equation and
system of\ infinity many\ parabolic equations is obtained. \ \ \ 

\begin{center}
\bigskip\ \ \textbf{AMS: }

\textbf{47Axx, 46E35, 47A50, 42B37, 42B15 }
\end{center}


\begin{center}
\textbf{1. Introduction }
\end{center}

Fourier multipliers in vector-valued function spaces have been studied e.g.
in $\left[ 4\right] $, $\left[ 18\right] ,$ $\left[ 31\right] .$
Operator-valued Fourier multipliers have been investigated in $\left[ 5%
\right] ,$ $\left[ 8-11\right] $ and $\left[ 29\right] .$ Mikhlin type
Fourier multipliers in scalar weighted spaces have been studied e.g. in $%
\left[ 13\right] $, $\left[ 28\right] $. Moreover, operator-valued Fourier
multipliers in weighted abstract $L_{p}$ spaces were investigated e.g. in $%
\left[ 2\right] $ and $\left[ 16\right] .$ In $\left[ \text{6, 12, 13}\right]
$ singular integral operators with operator-valued kernel were studied in
weighted $L_{p}$-spaces. Embedding theorems in vector-valued function spaces
are studied e.g. in $\left[ \text{14, 15}\right] $, $\left[ \text{20-26}%
\right] $. Regularity properties of differential-operator equations (DOEs)
have been studied e.g. in $\left[ 1\right] $, $\left[ 2\right] ,$ $\left[ 7,%
\text{ }8\right] ,$ $\left[ 22-25\right] $, $\left[ 29-30\right] .$ A
comprehensive introduction to DOEs and historical references may be found in 
$\left[ 1\right] $ and $\left[ 30\right] .$

In this paper, operator-valued multiplier theorems in $E-$valued weighted
Lebesque spaces are obtained. These multiplier theorems are used to show the
boundedness of embedding operator in the anisotropic Sobolev-Lions space $%
W_{p,\gamma }^{l}\left( \Omega ;E\left( A\right) ,E\right) $, i.e. under
some conditions we prove that the differential operator $u\rightarrow
D^{\alpha }u$ is bounded from $W_{p,\gamma }^{l}\left( \Omega ;E\left(
A\right) ,E\right) $ to $L_{p,\gamma }\left( \Omega ;E\left( A^{1-\left\vert
\alpha :l\right\vert }\right) \right) $ and the following
Ehrling-Nirenberg-Gagilardo type sharp estimate holds%
\begin{equation*}
\left\Vert D^{\alpha }u\right\Vert _{L_{p\mathbf{,\gamma }}\left( \Omega
;E\left( A^{1-\left\vert \alpha :l\right\vert -\mu }\right) \right) }\leq
C_{\mu }\left[ h^{\mu }\left\Vert u\right\Vert _{W_{p,\gamma }^{l}\left(
\Omega ;E\left( A\right) ,E\right) }+h^{-\left( 1-\mu \right) }\left\Vert
u\right\Vert _{L_{p\mathbf{,\gamma }}\left( \Omega ;E\right) }\right]
\end{equation*}%
for $u\in W_{p,\gamma }^{l}\left( \Omega ;E\left( A\right) ,E\right) $,
where $A$ is a positive operator in $E$ and%
\begin{equation*}
l=\left( l_{1},l_{2},...,l_{n}\right) ,\ \alpha =\left( \alpha _{1},\alpha
_{2},...,\alpha _{n}\right) ,\text{ }\left\vert \alpha :l\right\vert
=\dsum\limits_{k=1}^{n}\frac{\alpha _{k}}{l_{k}},
\end{equation*}%
\begin{equation*}
0<h\leq h_{0}<\infty ,\text{ }0<\mu <1-\left\vert \alpha :l\right\vert .
\end{equation*}%
This fact generalizes and improves the results $\left[ \text{3, \S\ 9, 27, 
\S\ 1.7}\right] $ for scalar Sobolev space, the result $\left[ 14\right] $
for one dimensional Sobolev-Lions spaces and the results $\left[ 15\right] $%
, $\left[ 22\right] $ for Hilbert-space valued class. Finally, we consider
the differential-operator equation%
\begin{equation}
\ Lu=\sum\limits_{\left\vert \alpha \right\vert =2l}a_{\alpha }D^{\alpha
}u+Au+\sum\limits_{\left\vert \alpha \right\vert <2l}A_{\alpha }\left(
x\right) D^{\alpha }u+\lambda u=f,  \tag{1.1}
\end{equation}%
where $a_{\alpha }$ are complex numbers, $A$, $A_{\alpha }\left( x\right) $
are linear operators in a Banach space $E$ and $\lambda $ is a complex
parameter.

We say that the problem $\left( 1.1\right) $ is $L_{p,\gamma }\left(
R^{n};E\right) $-separable if there exists a unique solution $u\in
W_{p,\gamma }^{2l}\left( R^{n};E\left( A\right) ,E\right) $ of $\left(
1.1\right) $ for all $f\in L_{p,\gamma }\left( R^{n};E\right) $ and there
exists a positive constant $C$ depend only on $p$ and $\gamma $ such that
the following coercive uniform estimate holds 
\begin{equation}
\sum\limits_{\left\vert \alpha \right\vert \leq 2l}\left\vert \lambda
\right\vert ^{1-\frac{\left\vert \alpha \right\vert }{2l}}\left\Vert
D^{\alpha }u\right\Vert _{L_{p,\gamma }\left( R^{n};E\right) }+\left\Vert
Au\right\Vert _{L_{p,\gamma }\left( R^{n};E\right) }\leq C\left\Vert
f\right\Vert _{L_{p,\gamma }\left( R^{n};E\right) }.  \tag{1.2}
\end{equation}

Estimate $\left( 1.2\right) $ implies that if $f\in L_{p,\gamma }\left(
R^{n};E\right) $ and $u$ is a solution of $\left( 1.1\right) $ then all
terms of equation $\left( 1.1\right) $ belong to $L_{p,\gamma }\left(
R^{n};E\right) $ (i.e. all terms are separable in $L_{p,\gamma }\left(
R^{n};E\right) $). The above estimate implies that the inverse of the
differential operator generated by $\left( 1.1\right) $ is bounded from $%
L_{p,\gamma }\left( R^{n};E\right) $ to $W_{p,\gamma }^{2l}\left(
R^{n};E\left( A\right) ,E\right) .$

By using the separability properties of $\left( 1.1\right) $ we show that
the Cauchy problem for the parabolic equation%
\begin{equation}
\ \partial _{t}u+\sum\limits_{\left\vert \alpha \right\vert =2l}a_{\alpha
}D^{\alpha }u+Au=f\left( t,x\right) ,\text{ }t\in \left( 0,\infty \right) 
\text{, }x\in R^{n},  \tag{1.3}
\end{equation}%
\begin{equation*}
u\left( 0,x\right) =0,\text{ }x\in R^{n}
\end{equation*}%
is well-posed in weighted spaces $L_{\mathbf{p},\gamma }\left(
R^{n};E\right) $ with mixed norm, where $\mathbf{p=}\left( p,p_{1}\right) $.

The paper is organized as follows. In Section 2, the necessary tools from
Banach space theory and some background materials are given. In sections 3,
the multiplier theorems in vector-valued weighted Lebesque spaces are
proved. In Section 4, by using these multiplier theorems, embedding theorems
in $E$-valued weighted Sobolev type spaces are shown. Finally, in sections
5-8 the separability properties of $\left( 1.1\right) $, $\left( 1.3\right) $
and also regularity properties of appropriate degenerate differential
operators are established.

\begin{center}
\textbf{2. Notations and background }
\end{center}

\ \ \ Let $E$ be a Banach space and let $\gamma =\gamma \left( x\right) ,$ $%
x=\left( x_{1},x_{2},...,x_{n}\right) $ be a positive measurable function on
the measurable subset $\Omega \subset R^{n}.$ Let $L_{p,\gamma }\left(
\Omega ;E\right) $ denote the weighted Lebesgue-Bochner space, i.e. the
space of all strongly measurable $E-$valued functions that are defined on $%
\Omega $ with the norm

\begin{equation*}
\left\Vert f\right\Vert _{L_{p,\gamma }}=\left\Vert f\right\Vert
_{L_{p,\gamma }\left( \Omega ;E\right) }=\left( \int \left\Vert f\left(
x\right) \right\Vert _{E}^{p}\gamma \left( x\right) dx\right) ^{\frac{1}{p}},%
\text{ }1\leq p<\infty ,
\end{equation*}

\begin{equation*}
\left\Vert f\right\Vert _{L_{\infty ,\gamma }\left( \Omega ;E\right) }=\text{%
ess}\sup\limits_{x\in \Omega }\left\Vert f\left( x\right) \right\Vert
_{E}\gamma \left( x\right) \text{ for }p=\infty .
\end{equation*}

For $\gamma \left( x\right) \equiv 1,$ the space $L_{p,\gamma }\left( \Omega
;E\right) $ will be denoted by $L_{p}=L_{p}\left( \Omega ;E\right) .$

The weight $\gamma $ is said to be satisfy an $A_{p}$ condition, i.e. $\
\gamma \in A_{p},$ $1<p<\infty $, if there is a positive constant $C$ such
that 
\begin{equation*}
\sup\limits_{Q}\left( \frac{1}{\left\vert Q\right\vert }\int\limits_{Q}%
\gamma \left( x\right) dx\right) \left( \frac{1}{\left\vert Q\right\vert }%
\int\limits_{Q}\gamma ^{-\frac{1}{p-1}}\left( x\right) dx\right) ^{p-1}\leq C
\end{equation*}%
for all for all cubes $Q$ $\subset R^{n}.$

The Banach space\ $E$\ is called a UMD-space and written as $E\in $ UMD if
only if the Hilbert operator 
\begin{equation*}
\left( Hf\right) \left( x\right) =\lim\limits_{\varepsilon \rightarrow
0}\int\limits_{\left\vert x-y\right\vert >\varepsilon }\frac{f\left(
y\right) }{x-y}dy
\end{equation*}%
is bounded in the space $L_{p}\left( R,E\right) ,$ $p\in \left( 1,\infty
\right) $ (see e.g. $\left[ 8\right] $). UMD spaces include $L_{p}$, $l_{p}$
spaces, Lorentz spaces $L_{pq},$ $p,$ $q\in \left( 1,\infty \right) $ and
Morrey spaces (see e.g. $\left[ 20\right] $).

A Banach space $E$ has a property ($\alpha $) (see e.g. $\left[ \text{19}%
\right] $) if there exists a constant $\alpha $ such that 
\begin{equation*}
\left\Vert \sum\limits_{i,j=1}^{N}\alpha _{ij}\varepsilon _{i}\varepsilon
_{j}^{\shortmid }x_{ij}\right\Vert _{L_{2}\left( \Omega \times \Omega
^{\shortmid };E\right) }dy\leq \alpha \left\Vert
\sum\limits_{i,j=1}^{N}\varepsilon _{i}\varepsilon _{j}^{\shortmid
}x_{ij}\right\Vert _{L_{2}\left( \Omega \times \Omega ^{\shortmid };E\right)
}
\end{equation*}%
for all $N\in \mathbb{N}\mathbf{,}$ $x_{i,j}\in E,$ $\alpha _{ij}\in \left\{
0,1\right\} ,$ $i,j=1,2,...,N,$ and all choices of independent, symmetric, $%
\left\{ -1,1\right\} -$ valued random variables $\varepsilon
_{1},\varepsilon _{2},...,\varepsilon _{N},$ $\varepsilon _{1}^{\prime
},\varepsilon _{2}^{\prime },...,\varepsilon _{N}^{\prime }$ \ on
probability spaces $\Omega ,$ $\Omega ^{\prime }.$ For example the spaces $%
L_{p}\left( \Omega \right) ,$ $1\leq p<\infty $ \ has the property ($\alpha $%
)$.$

Let $\mathbb{C}$ be the set of complex numbers and\ 

\begin{equation*}
\ S_{\varphi }=\left\{ \xi ;\text{ \ }\xi \in \mathbb{C}\text{, \ }%
\left\vert \arg \xi \right\vert \leq \varphi \right\} \cup \left\{ 0\right\}
,\text{ }0\leq \varphi <\pi .
\end{equation*}%
Let $E_{1}$ and\ $E_{2}$ be two Banach spaces. $B\left( E_{1},E_{2}\right) $
denotes the space of bounded linear operators from $E_{1}$ to $E_{2}.$ For $%
E_{1}=E_{2}=E$ it will be denote by $B\left( E\right) .$

A linear operator\ $A$ is said to be positive in\ a Banach\ space $E$,\ with
bound $M$, if\ $D\left( A\right) $ is dense in $E$ and 
\begin{equation*}
\left\Vert \left( A+\xi I\right) ^{-1}\right\Vert _{B\left( E\right) }\leq
M\left( 1+\left\vert \xi \right\vert \right) ^{-1}
\end{equation*}%
with $\xi \in S_{\varphi },$ $\varphi \in \left[ 0,\right. \left. \pi
\right) ,$ where $M$ is a positive constant and $I$ is an identity operator
in $E.$ Sometimes instead of $A+\xi I$,\ we will write $A+\xi $ or $A_{\xi
}. $ It is known $\left[ \text{27, \S 1.15.1}\right] $ there exist
fractional powers\ $A^{\theta }$ of the positive operator $A.$

\textbf{Definition 2.1}. A positive operator $A$ is said to be $R-$positive
in the Banach space $E$ if there exists $\varphi \in \left[ 0,\pi \right.
\left. {}\right) $ such that the set 
\begin{equation*}
\left\{ \xi \left( A+\xi I\right) ^{-1}:\xi \in S_{\varphi }\right\}
\end{equation*}%
is $R$-bounded (see e.g. $\left[ 29\right] $).

Let $E\left( A^{\theta }\right) $ denote the space $D\left( A^{\theta
}\right) $ with graphical norm defined as 
\begin{equation*}
\left\Vert u\right\Vert _{E\left( A^{\theta }\right) }=\left( \left\Vert
u\right\Vert ^{p}+\left\Vert A^{\theta }u\right\Vert ^{p}\right) ^{\frac{1}{p%
}},\text{ }1\leq p<\infty ,\text{ }-\infty <\theta <\infty .
\end{equation*}%
Let $\left( E_{1},E_{2}\right) _{\theta ,p}$denote the interpolation space\
obtained from $\left\{ E_{1},E_{2}\right\} $ by the $K-$method $\left[ \text{%
27, \S 1.3.1}\right] $, where $\theta \in \left( 0,1\right) ,$ $p\in \left[
1,\infty \right. \left. {}\right) $.

We denote by $D\left( R^{n};E\right) $ the space of $E-$valued $C^{\infty }-$
function with compact support, equipped with the usual inductive limit
topology and $S\left( E\right) =S\left( R^{n};E\right) $ denote the $E-$%
valued Schwartz space of rapidly decreasing smooth functions. For $E=\mathbb{%
C}$ we simply write $D\left( R^{n}\right) $ and $S=S\left( R^{n}\right) $,
respectively. Let $D^{\prime }\left( R^{n};E\right) $ $=B\left( D\left(
R^{n}\right) ,E\right) $ denote the space of $E-$valued distributions and
let $S^{\prime }\left( E\right) =S^{\prime }\left( R^{n};E\right) $ denote a
space of linear continued mapping from $S\left( R^{n}\right) $ into\ $E.$
The Fourier transform for $u\in S^{\prime }\left( R^{n};E\right) $ is
defined by 
\begin{equation*}
F\left( u\right) \left( \varphi \right) =u\left( F\left( \varphi \right)
\right) \text{, }\varphi \in S\left( R^{n}\right) .
\end{equation*}%
Let $\gamma $ be such that $S\left( R^{n};E_{1}\right) $ is dense in $%
L_{p,\gamma }\left( R^{n};E_{1}\right) .$ A function 
\begin{equation*}
\Psi \in C^{\left( l\right) }\left( R^{n};B\left( E_{1},E_{2}\right) \right)
\end{equation*}%
is called a multiplier from\ $L_{p,\gamma }\left( R^{n};E_{1}\right) $ to $%
L_{q,\gamma }\left( R^{n};E_{2}\right) $ if there exists a positive constant 
$C$ such that 
\begin{equation*}
\left\Vert F^{-1}\Psi \left( \xi \right) Fu\right\Vert _{L_{q,\gamma }\left(
R^{n};E_{2}\right) }\leq C\left\Vert u\right\Vert _{L_{p,\gamma }\left(
R^{n};E_{1}\right) }
\end{equation*}%
for all $u\in S\left( R^{n};E_{1}\right) $.

We denote the set of all multipliers fom\ $L_{p,\gamma }\left(
R^{n};E_{1}\right) $ to $L_{q,\gamma }\left( R^{n};E_{2}\right) $ by $%
M_{p,\gamma }^{q,\gamma }\left( E_{1},E_{2}\right) .$ For $E_{1}=E_{2}=E$ we
denote the $M_{p,\gamma }^{q,\gamma }\left( E_{1},E_{2}\right) $ by $%
M_{p,\gamma }^{q,\gamma }\left( E\right) .$

A set $K\subset B(E_{1},E_{2})$ is called $R$-bounded (see e.g. $[$8, \S\ 3.1%
$]$) if there is a constant $C>0$ such that for all $T_{1},T_{2},...,T_{m}%
\in K$ and $u_{1},u_{2},...,u_{m}\in E_{1},$ $m\in \mathbb{N}$

\begin{equation*}
\int\limits_{0}^{1}\left\Vert
\sum\limits_{j=1}^{m}r_{j}(y)T_{j}u_{j}\right\Vert _{E_{2}}dy\leq
C\int\limits_{0}^{1}\left\Vert \sum\limits_{j=1}^{m}r_{j}(y)u_{j}\right\Vert
_{E_{1}}dy,
\end{equation*}%
where $\left\{ r_{j}\right\} $ is a sequence of independent symmetric $%
\left\{ -1;1\right\} $-valued random variables on$\left[ 0,1\right] .$ The
smallest $C$ for which the above estimate holds is called the $R-$bound of $%
K $ and denoted by $R\left( K\right) .$

\textbf{Definition 2.2. }The Banach space $E$ satisfies the multiplier
condition with respect to\ $p\in \left( 1,\infty \right) $ and to the
weighted function $\gamma $\ if for all $\Psi \in C^{\left( n\right) }\left(
R^{n};B\left( E\right) \right) $ the inequality

\begin{equation}
R\left\{ \left\vert \xi \right\vert ^{\left\vert \alpha \right\vert }D_{\xi
}^{\alpha }\Psi \left( \xi \right) :\xi \in R^{n}\smallsetminus \left\{
0\right\} ,\right\} \leq K_{\alpha }<\infty  \tag{2.1}
\end{equation}%
for $\alpha =\left( \alpha _{1},\alpha _{2},...,\alpha _{n}\right) ,$ $%
\alpha _{i}\in \left\{ 0,1\right\} $ implies that $\Psi \in M_{p,\gamma
}^{p,\gamma }\left( E\right) .$

Note that, if $E_{1}$ and $E_{2}$ are UMD spaces and $\gamma \left( x\right)
\equiv 1,$ then by virtue of operator valued multiplier theorems (see e.g $%
\left[ 9-12\right] ,$ $\left[ 30\right] $) we obtain that $\Psi $ is a
Fourier multiplier in $L_{p}\left( R^{n};E\right) .$

Let$\ \Omega $ be a domain on $R^{n}$ and let $l=\left(
l_{1},l_{2},...,l_{n}\right) \in \mathbb{N}^{n}$. Assume $E_{0}$ is
continuously and densely belongs to $E.$ Here, $W_{p,\gamma }^{l}\left(
\Omega ;E_{0},E\right) $ denotes the anisotropic weighted Sobolev-Lions type
space of functions $u\in L_{p,\gamma }\left( \Omega ;E_{0}\right) $ which
have generalized derivatives $\frac{\partial ^{l_{k}}u}{\partial
x_{k}^{l_{k}}}\in L_{p,\gamma }\left( \Omega ;E\right) $ with norm 
\begin{equation*}
\left\Vert u\right\Vert _{W_{p,\gamma }^{m}\left( \Omega ;E_{0},E\right)
}=\left\Vert u\right\Vert _{L_{p,\gamma }\left( \Omega ;E_{0}\right)
}+\sum\limits_{k=1}^{n}\left\Vert \frac{\partial ^{l_{k}}u}{\partial
x_{k}^{l_{k}}}\right\Vert _{L_{p,\gamma }\left( \Omega ;E\right) }<\infty .
\end{equation*}

\ \ \ \ For $l_{1}=l_{2}=...=l_{n}=m$ we denote $W_{p,\gamma }^{l}\left(
\Omega ;E_{0},E\right) $ by $W_{p,\gamma }^{m}\left( \Omega ;E_{0},E\right) $
as a isotropic weighted Sobolev-Lions space.

\bigskip

\begin{center}
\textbf{3. Operator-valued multiplier results in weighted Lebesque spaces}
\end{center}

\bigskip Let $E_{1}$, $E_{2}$ be\ Banach spaces. We put 
\begin{equation*}
X=L_{p,\gamma }\left( R^{n};E_{1}\right) \text{ and }Y=L_{p,\gamma }\left(
R^{n};E_{2}\right) .
\end{equation*}

By following Theorems 3. 6 and 3.7 of $\left[ 9\right] $ we will show the
following multiplier theorems:

\textbf{Theorem 3.1}$.$ Let $\gamma \in A_{p},$ $p\in \left( 1,\infty
\right) $. Assume $E_{1}$, $E_{2}$ are\ UMD spaces with property ($\alpha $)
and let 
\begin{equation*}
M\in C^{\left( n\right) }\left( R^{n}\smallsetminus \left\{ 0\right\}
;B\left( E_{1},E_{2}\right) \right) .
\end{equation*}%
If 
\begin{equation*}
R\left\{ \xi ^{\beta }D_{\xi }^{\beta }M\left( \xi \right) :\xi \in
R^{n}\smallsetminus \left\{ 0\right\} \right\} \leq C_{\beta }<\infty
\end{equation*}%
for all $\beta =\left( \beta _{1},\beta _{2},...,\beta _{n}\right) ,$ $\beta
_{i}\in \left\{ 0,1\right\} ,$ then $M$ is a multiplier from $X$ to $Y$ with 
$\left\Vert M\right\Vert _{B\left( X,Y\right) }\leq C\dsum\limits_{\beta
_{i}\in \left\{ 0,1\right\} }C_{\beta }.$

If $n=1$ the result remains true without $E_{1}$ having property ($\alpha $).

\textbf{Theorem 3.2}$.$ Let $\gamma \in A_{p},$ $p\in \left( 1,\infty
\right) $. Let $E_{1}$, $E_{2}$ be\ UMD spaces and let 
\begin{equation*}
M\in C^{\left( n\right) }\left( R^{n}\smallsetminus \left\{ 0\right\}
;B\left( E_{1},E_{2}\right) \right) .
\end{equation*}%
If 
\begin{equation*}
R\left\{ \left\vert \xi \right\vert ^{\left\vert \beta \right\vert }D_{\xi
}^{\beta }M\left( \xi \right) :\xi \in R^{n}\smallsetminus \left\{ 0\right\}
\right\} \leq C_{\beta }<\infty
\end{equation*}%
for all $\beta =\left( \beta _{1},\beta _{2},...,\beta _{n}\right) ,$ $\beta
_{i}\in \left\{ 0,1\right\} ,$ then $M$ is a\ multiplier from $X$ to $Y$\
with $\left\Vert M\right\Vert _{B\left( X,Y\right) }\leq
C\dsum\limits_{\beta _{i}\in \left\{ 0,1\right\} }C_{\beta }.$

\bigskip To prove Theorem 3.1 we need the following result:

The following Propositions A$_{1}$ and A$_{2}$ are due to Cl\'{e}ment, de
Pagter, Sukochev and Witvliet, see $[5]$.

\textbf{Proposition A}$_{1}$\textbf{. }\ Let $\Delta _{j}^{E_{1}}$ and $%
\Delta _{j}^{E_{2}}$ be unconditional Schauder decompositions of the Banach
spaces $E_{1}$ and $E_{2}$ respectively, with unconditional constants $%
C_{E_{1}}$ and $C_{E_{2}}$. Further let $\left\{ T_{j}:j\in \mathbb{Z}%
^{n}\right\} $ be an $R-$bounded family in $B(E_{1},E_{2})$ with $%
T_{j}\Delta _{j}^{E_{1}}=\Delta _{j}^{E_{2}}T_{j}\Delta _{j}^{E_{1}}$ for
all $j\in \mathbb{N}.$ Then the series 
\begin{equation*}
Tu=\dsum\limits_{j=1}^{\infty }T_{j}\Delta _{j}^{E_{1}}u
\end{equation*}%
converges for every $u\in E_{1}$ and defines a bounded operator $%
T:E_{1}\rightarrow E_{2}$ with 
\begin{equation*}
\left\Vert T\right\Vert \leq C_{E_{1}}C_{E_{2}}R\left( \left\{ T_{j}:j\in 
\mathbb{Z}^{n}\right\} \right) .
\end{equation*}

\textbf{Proposition A}$_{2}$\textbf{. }Assume $E$ is a Banach space that has
property($\alpha $), $\Delta =\left\{ \Delta _{k}\right\} _{k=1}^{\infty }$
is an unconditional Schauder decomposition and $Q\subset B\left( E\right) $
is an $R$-bounded collection of operators. Then the set 
\begin{equation*}
S:=\left\{ \dsum\limits_{k=0}^{\infty }T_{k}\Delta _{k}:T_{k}\in Q\text{
such that }T_{k}\Delta _{k}=\Delta _{k}T_{k}\text{ for all }k\in \mathbb{N}%
\right\}
\end{equation*}%
is $R$-bounded in $E$.

Let $\Omega \subset R^{n}$. By using the same reasoning as used in $\left[ 
\text{5, Lemma 3.17}\right] $ we have:

\textbf{Lemma 3.1}. Let $\gamma \in A_{p}$, $p\in \left( 1,\infty \right) $.
Assume $E$ is a Banach spaces. For $\phi \in L_{\infty }\left( \Omega
\right) $ we denote by $M_{\phi }=M_{\phi }^{X}$ the associated
multiplication operator in $X$ $=L_{p,\gamma }\left( \Omega ;E\right) $.
Then the collection 
\begin{equation*}
\left\{ M_{\phi }:\phi \in L_{\infty }\left( \Omega \right) \text{, }%
\left\Vert \phi \right\Vert _{\infty }\leq 1\right\}
\end{equation*}%
is $R-$bounded in $X.$

From Lemma 3.1 we obtain

\textbf{Corollary 3.1. }Let $\gamma \in A_{p},$ $p\in \left( 1,\infty
\right) $.\textbf{\ }Assume $E_{1}$ and $E_{2}$ are Banach spaces. For $\phi
\in L_{\infty }\left( \Omega \right) $ we denote by $M_{\phi }^{X}$ and $%
M_{\phi }^{Y}$ the associated multiplication operators in $X=L_{p,\gamma
}\left( \Omega ;E_{1}\right) $ and $Y=L_{p,\gamma }\left( \Omega
;E_{2}\right) $ respectively. If the set $K\subset B(X,Y)$ is $R-$bounded,
then the family 
\begin{equation*}
\left\{ M_{\phi }^{X}TM_{\phi }^{Y}:\phi \text{, }\psi \in L_{\infty }\left(
R^{n}\right) \text{, }\left\Vert \phi \right\Vert _{\infty },\text{ }%
\left\Vert \psi \right\Vert _{\infty }\leq 1\text{, }T\in K\right\}
\end{equation*}%
is $R-$bounded as well.

\bigskip For $k=nr+j$, $r\in \mathbb{Z}$, $j\in \left\{ 1,2,...,n\right\} $
let 
\begin{equation*}
\mathbb{D}_{k}=\left\{ \xi =\left( \xi _{1},\xi _{2},...,\xi _{n}\right) \in
R^{n}\right. ,\text{ }\left\vert \xi _{i}\right\vert <2^{r+1}\text{ for }%
i\in \left\{ 1,2,...,j-1\right\} ,
\end{equation*}%
\begin{equation*}
2^{r}\leq \left\vert \xi _{j}\right\vert <2^{r+1},\text{ }\left\vert \xi
_{i}\right\vert <2^{r}\text{ \ for}\left. i\in \left\{ j+1,...,n\right\}
\right\} .
\end{equation*}%
For $\nu =\left( \nu _{1},\nu _{2},...,\nu _{n}\right) \in \mathbb{Z}^{n}$
let%
\begin{equation*}
\mathbb{\Delta }_{\nu }=\left\{ \xi \in R^{n}\smallsetminus \left\{
0\right\} \text{, }2^{\nu _{i}-1}\leq \left\vert \xi _{j}\right\vert <2^{\nu
_{i}}\text{ for }i\in \left\{ 1,2,...,n\right\} \right\} .
\end{equation*}

From $\left[ \text{2, Proposition A}_{4}\right] $ we have:

\textbf{Lemma 3.2}. Let $\gamma \in A_{p}$, $p\in \left( 1,\infty \right) $\
and let $E$ be a UMD space (respectively, UMD space with property ($\alpha $%
)). Then for any choice of signs $\varepsilon _{k}$, $k\in \mathbb{Z}$
(respectively, $\varepsilon _{k}$, $k\in \mathbb{Z}^{n}$ ) the function $%
\psi :R^{n}\rightarrow \mathbb{C}$ with $\psi \left( \xi \right)
=\varepsilon _{k}$ for $\xi \in \mathbb{D}_{k}$ (respectively, $\psi \left(
\xi \right) =\varepsilon _{\nu }$ for $\xi \in \mathbb{D}_{\nu },$ $\nu \in 
\mathbf{\bigtriangleup }_{\nu }$) is a $M_{p,\gamma }^{p,\gamma }\left(
E\right) $\ multiplier.

\bigskip Let $E$ be a Banach space. The ($n-$dimensional) Riesz projection
operator $R$ is defined by 
\begin{equation*}
Rf=F^{-1}\chi _{\left( 0,\infty \right) ^{n}}Ff\text{, \ }f\in S\left(
R^{n};E\right) ,
\end{equation*}%
where $\chi \left( \Omega \right) $ denotes the characteristic function of $%
\Omega $ $\subset R^{n}.$

Let 
\begin{equation*}
R_{j}f=F_{j}^{-1}\chi _{j}F_{j}f\text{ for }f\in S\left( R^{n};E\right) 
\text{, }j=1,2,...,n,
\end{equation*}%
where $F_{j}$ denote the one-dimensional Fourier transform with respect to
variable $x_{j}$ and $\chi _{j}$ denotes the characteristic function of the
halfspace 
\begin{equation*}
R_{j}^{n}=\left\{ x=\left( x_{1},x_{2},...,x_{n}\right) \in R^{n}\text{, }%
x_{j}>0\right\} .
\end{equation*}

\textbf{Lemma 3.3}. Assume $\gamma \in A_{p}$ for $p\in \left( 1,\infty
\right) $\ and $E$ is a UMD space. Then $R$ defines a bounded operator in $%
L_{p,\gamma }\left( R^{n};E\right) .$

\textbf{Proof. }Since $\gamma \in A_{p}$, then by $\left[ \text{12,
Corollary 2.10}\right] $ (or $\left[ \text{6, Theorem 4}\right] $ )\ the
Hilbert operator is bounded in $L_{p,\gamma }\left( R;E\right) .$ It is
known that $R_{1}=\frac{1}{2\pi i}\left( i\pi I-H\right) $, where $I$ is the
identity operator. By using this relation we obtain that Riesz projection
operator $R_{1}$ is bounded in $L_{p,\gamma }\left( R;E\right) $. Hence,
one-dimensional Riesz projection $R_{j}$ also are defined bounded operators
in $L_{p,\gamma }\left( R;E\right) $. It is not hard to see that 
\begin{equation*}
R=\dprod\limits_{j=1}^{n}R_{j},
\end{equation*}%
i.e. $R$ is bounded operator in $L_{p,\gamma }\left( R;E\right) .$

For $j=\left( j_{1},j_{2},...j_{n}\right) \in \mathbb{Z}^{n}$ let $D_{j}$ be
the dyadic interval associated with $j$, i.e.%
\begin{equation}
D_{j}=\prod\limits_{k=1}^{n}\left[ 2^{j_{k}},\left. 2^{j_{k}+1}\right)
\right.  \tag{3.1}
\end{equation}%
and 
\begin{equation*}
\text{ }Q=Q_{a,b}=\prod\limits_{k=1}^{n}\left( a_{k},b_{k}\right) ,
\end{equation*}%
where 
\begin{equation*}
\text{ }a=\left( a_{1},a_{2},...,a_{n}\right) \text{ and }b=\left(
b_{1},b_{2},...,b_{n}\right) \in R^{n}.
\end{equation*}%
Consider the operator 
\begin{equation*}
\Phi _{a,b}f=F^{-1}\chi \left( Q_{a,b}\right) Ff\text{ }\ \text{for }f\in
S\left( R^{n};E\right) .
\end{equation*}

\textbf{Lemma 3.4}. Assume $\gamma \in A_{p}$ for $p\in \left( 1,\infty
\right) $\ and $E$ is a UMD space. Then for each $a,b$ $\in R^{n}$ the
operator $f\rightarrow \Phi _{a,b}f$ is bounded in $L_{p,\gamma }\left(
R^{n};E\right) $. Moreover, the set $\left\{ \Phi _{a,b}:a,b\in
R^{n}\right\} $ is an $R-$bounded family in $B\left( L_{p,\gamma }\left(
R^{n};E\right) \right) $.

\textbf{Proof. }We first look at characteristic functions of sets of the form%
\begin{equation*}
C_{a}=\prod\limits_{k=1}^{n}\left[ a_{k},\right. \left. \infty \right) .
\end{equation*}%
We can $F^{-1}\chi _{C_{a}}Ff$ expressed as:%
\begin{equation*}
\Phi _{a}f=F^{-1}\chi _{C_{a}}Ff=e^{ia_{1}\tau _{1}}R_{1}e^{-ia_{1}\tau
_{1}}...e^{ia_{n}\tau _{n}}R_{1}e^{-ia_{n}\tau _{n}}f\left( \tau
_{1},...,\tau _{n}\right)
\end{equation*}%
for%
\begin{equation*}
\tau _{1},...,\tau _{n}\in R^{n}.
\end{equation*}%
We see that the set $\left\{ \Phi _{a}:a\in R^{n}\right\} $ is $R-$bounded
in view of Proposition 3.1. Setting $C_{b}=\prod\limits_{k=1}^{n}\left[
-\infty ,\right. \left. b_{k}\right) $ we analogously get that the set $%
\left\{ \Phi _{b}:b\in R^{n}\right\} $ is $R-$bounded as well, where 
\begin{equation*}
\Phi _{b}f=F^{-1}\chi _{C_{b}}Ff\text{ for }f\in S\left( R^{n};E\right) .
\end{equation*}%
Since $\Phi _{a,b}=\Phi _{a}\Phi _{b},$ the result follows because the
pointwise product of $R-$bounded sets is again $R-$bounded.

Assume $E_{1}$ and $E_{2}$ are UMD spaces. We put 
\begin{equation*}
X_{1}=R\left( L_{p,\gamma }\left( R^{n};E_{1}\right) \right) ,\text{ }%
Y_{1}=R\left( L_{p,\gamma }\left( R^{n};E_{2}\right) \right) .
\end{equation*}

\bigskip Let $\left\{ A_{j}:j\in \mathbb{Z}^{n}\right\} $\ be a
decomposition of $(0,\infty )^{n}$ in intervals such that for each compact $%
K\in (0,\infty )^{n}$ the set $\left\{ A_{j}\cap K:j\in \mathbb{Z}%
^{n}\right\} $ is finite. Assume further that the families $\left\{ \Delta
_{j}^{X_{1}}:j\in \mathbb{Z}^{n}\right\} $ and $\left\{ \Delta
_{j}^{Y_{1}}:j\in \mathbb{Z}^{n}\right\} $ of the corresponding Fourier
multipliers, i.e 
\begin{equation*}
\Delta _{j}^{X_{1}}=F_{E_{1}}\chi _{A_{j}}F_{E_{1}}^{-1}\text{, }\Delta
_{j}^{Y_{1}}=F_{E_{2}}\chi _{A_{j}}F_{E_{2}}^{-1}
\end{equation*}%
are unconditional Schauder decompositions of $X_{1}$, $Y_{1}$ respectively,
where $F_{E}$ and $F_{E}^{-1}$ denote the Fourier and inverse Fourier
transforms. For $k\in \mathbb{N}$ we now cut each interval $A_{j}$ in $%
2^{kn} $ smaller ones by decomposing it in each coordinate direction into $%
2^{k}$ pieces. These new smaller intervals are denoted by $A_{j,l}^{k},$
where $j\in \mathbb{Z}^{n}$ and $l\in \left\{ 0,1,...,2^{k}-1\right\} ^{n}.$

Let $M$ be a function on $R^{n}$ with values in a Banach space $%
B(E_{1},E_{2})$. Assume that $M$ is constant operator on the intervals $%
A_{j,l}^{k}$, and denote by $M_{j,l}^{k}$, the corresponding value of $M$.
Next we show that an operator-valued function which is constant on the $%
A_{j,l}^{k}$ 's is a Fourier multiplier from $X$ to $Y$ if it satisfies a
certain inequality involving $R-$bounds.

\textbf{Proposition 3.1. }Assume $\gamma \in A_{p}$ for $p\in \left(
1,\infty \right) $\ and $E_{1}$, $E_{2}$ are UMD spaces. Further let $M:$ $%
R^{n}$ $\rightarrow B(E_{1},E_{2})$ be a function which is constant on each $%
A_{j,l}^{k}$ and zero on $R^{n}\smallsetminus \left( 0,\infty \right) ^{n}.$
Assume that 
\begin{equation*}
\dsum\limits_{r=\alpha }^{\beta .\left( 2^{k}-1\right) }R\left( \left\{
\dsum\limits_{\nu \in \left( 0,1\right) ^{n},\nu \leq \beta }\left(
-1\right) ^{\left\vert \nu \right\vert }M_{j,\beta \left( r-\nu \right)
}^{k}:j\in \mathbb{Z}^{n}\right\} \right) =C_{\beta ,k}<\infty
\end{equation*}%
for every multiindex $\beta \in \left( 0,1\right) ^{n}$ and $k\in \mathbb{Z}%
. $ Then $M$ is a Fourier multiplier from $X$ into $Y$. The norm of $%
T=F_{E_{2}}^{-1}MF_{E_{1}}$ may be estimated by 
\begin{equation*}
\left\Vert T\right\Vert \leq C_{X}C_{Y}C_{Q}\dsum\limits_{\beta \in \left(
0,1\right) ^{n}}C_{\beta ,k}
\end{equation*}%
where $C_{X}$ and $C_{Y}$ are the unconditional constants and $C_{Q}$ is the 
$R$--bound found in Lemma 3.4.

\textbf{Proof. }By Lemma 3.4, each $\chi _{A_{j,l}^{k}}$\ is a Fourier
multiplier in $X$. We denote the operators $F_{E_{1}}\chi
_{A_{j,l}^{k}}F_{E_{1}}^{-1}$ by $\Delta _{j,l}^{k}.$ For $f\in S\left(
R^{n};E_{1}\right) $ we get 
\begin{equation*}
Tf=F_{E_{2}}^{-1}MF_{E_{1}}f=F_{E_{2}}^{-1}\dsum\limits_{j=-\infty }^{\infty
}M_{\chi _{A_{j}}}F_{E_{1}}f.
\end{equation*}

Then by using the same reasoning as used in $\left[ \text{9, Theorem 3.3}%
\right] $ we obtain 
\begin{equation*}
Tf=\dsum\limits_{j=-\infty }^{\infty }T_{j}\Delta _{j}^{X}f\text{,}
\end{equation*}%
where $T_{j}$ are operators defined by 
\begin{equation*}
T_{j}=\dsum\limits_{\beta \in \left( 0,1\right) ^{n}}\dsum\limits_{r=\alpha
}^{\beta .\left( 2^{k}-1\right) }\dsum\limits_{\nu \in \left( 0,1\right)
^{n},\nu \leq \beta }\left( -1\right) ^{\left\vert \nu \right\vert
}M_{j,\beta \left( r-\nu \right) }^{k}\dsum\limits_{l=r}^{2^{k}-1}\Delta
_{j,l}^{k}.
\end{equation*}

Since $M_{j,l}^{k}\Delta _{j,l}^{k}=\Delta _{j}^{Y}M_{j,l}^{k}$ and $\Delta
_{j,l}^{k}\Delta _{j}^{X}=\Delta _{j}^{X}\Delta _{j,l}^{k},$ we have $\Delta
_{j,l}^{k}\Delta _{j}^{X}=T_{j}\Delta _{j}^{X}=\Delta _{j}^{Y}T_{j}\Delta
_{j}^{X}.$ Moreover, since $\left\{ \Delta _{j}^{X}:j\in \mathbb{Z}%
^{n}\right\} $\ and $\left\{ \Delta _{j}^{Y}:j\in \mathbb{Z}^{n}\right\} $
are unconditional Schauder decompositions of the spaces $X$, $Y$
respectively and $S\left( R^{n};E_{1}\right) $ is dense in $X$, it remains
to prove that the family $\left\{ T_{j}:j\in \mathbb{Z}^{n}\right\} $ is $R-$%
bounded. This step is derived as in $\left[ \text{9, Theorem 3.3}\right] $,
i.e. we show that%
\begin{equation*}
R\left( \left\{ T_{j}:j\in \mathbb{Z}^{n}\right\} \right) \leq
C_{Q}\dsum\limits_{\beta \in \left( 0,1\right) ^{n}}C_{\beta ,k}.
\end{equation*}

Then in view of Proposition A$_{1}$ we have $T\in B\left( X;Y\right) $ with 
\begin{equation*}
\left\Vert T\right\Vert \leq C_{X}C_{E}R\left( \left\{ T_{j}:j\in \mathbb{Z}%
^{n}\right\} \right) \leq C_{X}C_{E}C_{Q}\dsum\limits_{\beta \in \left(
0,1\right) ^{n}}C_{\beta ,k}.
\end{equation*}

In a similar way as $\left[ \text{9, Proposition 3.4}\right] $ it can be
shown the following proposition. It will be used to prove the Mikhlin
theorem by approximating the given function $\Psi :R^{n}\rightarrow B(X,Y)$
by piecewise constant multipliers and is a generalization of the same result
from $\left[ 7\right] $ for unweighted spaces $L_{p}\left( R^{n};E\right) $.

\bigskip \textbf{Proposition 3.2. }Assume $\gamma \in A_{p}$ for $p\in
\left( 1,\infty \right) $\ and $E_{1}$, $E_{2}$ are Banach spaces. Let $M$, $%
M_{N}\in L_{1}^{loc}\left( R^{n},B\left( E_{1},E_{2}\right) \right) $ be
Fourier multipliers from $X$ to $Y$ such that $M_{N}\rightarrow M$ in $%
L_{1}^{loc}\left( R^{n},B\left( E_{1},E_{2}\right) \right) $. If $E_{2}$
reflexive and\ the sequence 
\begin{equation*}
\left\{ T_{N}\right\} =\left\{ F^{-1}M_{N}F\text{, }N\in \mathbb{N}\right\}
\end{equation*}%
is uniformly bounded in $B\left( X,Y\right) ,$ then the operator $%
T:=F_{E_{2}}^{-1}MF_{E_{1}}$ is a bounded operator from $X$ to $Y$ with 
\begin{equation*}
\left\Vert T\right\Vert \leq \lim\limits_{N\rightarrow \infty }\text{inf}%
\left\Vert T_{N}\right\Vert .
\end{equation*}

The next lemma states that the family of dyadic intervals in $R^{n}$ can be
used to build up an unconditional Schauder decomposition of $R\left(
X\right) $ provided $E$ is a UMD space with property ($\alpha $).

\textbf{Lemma 3.5}. Assume $\gamma \in A_{p}$ for $p\in \left( 1,\infty
\right) $\ and $E$ is a UMD space. For $j=\left( j_{1},j_{2},...j_{n}\right)
\in \mathbb{Z}^{n}$ let $D_{j}$ be the dyadic interval defined by $\left(
3.1\right) $ and $\Delta _{j}:=F^{-1}\chi _{D_{j}}F.$ Then:

(a) If $n=1$, then the family $\left\{ \Delta _{j}:j\in \mathbb{Z}%
^{n}\right\} $ is an unconditional Schauder decomposition of $X_{1}=R\left(
L_{p,\gamma }\left( R^{n};E\right) \right) ;$

(b) If $E$ has property ($\alpha $), then the assertion of part (a) is true
for arbitrary $n$.

\textbf{Proof. }(a) It is clear that the $\Delta _{j}$'s are projections in $%
L_{p,\gamma }\left( R^{n};E\right) $ and that $\Delta _{j}\Delta _{j^{\prime
}}$ $=$ $\delta _{jj^{\prime }}\Delta _{j}$. Let $1,2,...$ be any
enumeration of $\mathbb{Z}$. We have to prove that 
\begin{equation*}
T_{N}f:=\dsum\limits_{k=1}^{N}\Delta _{l_{k}}f\rightarrow f\text{ \ in }X_{1}%
\text{ as }N\rightarrow \infty .
\end{equation*}

This convergence is clear for $f\in $ $S(0,\infty ;E)$. In view of a $%
3\varepsilon $--argument it remains to show that the set$\{T_{N}:N\in 
\mathbb{N}\}$ is uniformly bounded. To this aim we define the function $%
m_{N}:\mathbb{R}\rightarrow \mathbb{R}$ by 
\begin{equation*}
m_{N}\left( x\right) =\left\{ 
\begin{array}{c}
1\text{ when }x\in \cup _{k=1}^{N}D_{l_{k}}, \\ 
-1\text{ when }x\text{ is elsewhere and }N\in \mathbb{N}%
\end{array}%
\right. .
\end{equation*}

By Proposition A$_{4}$ of $[2]$ we get that each $m_{N}\left( x\right) $ is
a Fourier multiplier in $L_{p,\gamma }\left( R^{n};E\right) .$ Moreover, the
proof the Proposition A$_{4}$ in $[2]$ shows that the family $\left\{
F^{-1}m_{N}F\right\} $ is uniformly bounded. Hence, we get 
\begin{equation*}
\left\Vert T_{N}\right\Vert =\left\Vert \dsum\limits_{k=1}^{N}F^{-1}\chi
_{D_{l_{k}}}F\right\Vert =\left\Vert \dsum\limits_{k=1}^{N}F^{-1}\chi _{\cup
_{k=1}^{N}D_{l_{k}}}F\right\Vert =
\end{equation*}%
\begin{equation*}
\frac{1}{2}\left\Vert \dsum\limits_{k=1}^{N}F^{-1}\left( \chi _{\left(
0,\infty \right) }+m_{N}\right) F\right\Vert \leq \frac{1}{2}\left(
\left\Vert R\right\Vert +\sup\limits_{N\in \mathbb{N}}\left\Vert
F^{-1}m_{N}F\right\Vert \right) <\infty .
\end{equation*}

This gives the assertion (a). By Proposition A$_{2}$ we get that the
collection $\left\{ \dsum\limits_{k\in G}\Delta _{k}:G\subset \mathbb{Z}%
\right\} $ is $R-$bounded which in view of Proposition A$_{1}$ yields that
the product of two unconditional Schauder decompositions is again an
unconditional Schauder decomposition. The general case now follows by
induction.

\textbf{Proof of Theorem 3.1. }Without loss of generality we assume $M\left(
\xi \right) =0$ for $\xi \notin \left( 0,\infty \right) ^{n}.$ To apply
Propositions 3.1 and 3.2 we use the decomposition of Lemma 3.5 to
approximate $M$. Now, we cut each $D_{j}$ into $2^{nk}$ pieces and define 
\begin{equation*}
M_{j,r}^{k}:=M\left(
2^{j_{1}}+r_{1}2^{j_{1}-k},...,2^{j_{n}}+r_{n}2^{j_{n}-k}\right) \text{, }%
k\in \mathbb{Z}\text{, }r\text{, }j\in \mathbb{Z}^{n},
\end{equation*}%
where 
\begin{equation*}
0\leq r_{i}\leq 2^{k}-1.
\end{equation*}

In view of Proposition 3.1 we have to estimate the $R-$bounds%
\begin{equation*}
\dsum\limits_{r=\alpha }^{\beta .\left( 2^{k}-1\right) }R\left( \left\{
\dsum\limits_{\nu \in \left( 0,1\right) ^{n},\nu \leq \beta }\left(
-1\right) ^{\left\vert \nu \right\vert }M_{j,\beta \left( r-\nu \right)
}^{k}:j\in \mathbb{Z}^{n}\right\} \right)
\end{equation*}%
for all $\beta \in \left( 0,1\right) ^{n}$ independently of $k$. For $\beta
= $ $\left( 0,0,...,0\right) $ this expression is trivially bounded by $%
R(\{M(\xi ),$ $\xi \neq 0\})$. For $\beta \neq 0$ let $i$ be the smallest
index with $\beta _{i}$ $=1$. Every $\nu $ with $\nu _{i}=0$ and $\nu \leq
\beta $ has a term $\tilde{\nu}$ with $\tilde{\nu}_{m}=\nu _{m}$ for $m\neq
i $ and $\tilde{\nu}_{i}=1$. Now, by using the same reasoning as used in the
proof of Theorem 3.6 of $\left[ 9\right] $ by Corollary 3.1 we get the
desired estimate 
\begin{equation*}
\dsum\limits_{r=\alpha }^{\beta .\left( 2^{k}-1\right) }R\left( \left\{
\dsum\limits_{\nu \in \left( 0,1\right) ^{n},\nu \leq \beta }\left(
-1\right) ^{\left\vert \nu \right\vert }M_{j,\beta \left( r-\nu \right)
}^{k}:j\in \mathbb{Z}^{n}\right\} \right) \leq
\end{equation*}%
\begin{equation*}
CR\left( \left\{ \xi ^{\beta }D^{\beta }M:\xi \in \left( 0,\infty \right)
^{n}\right\} \right) \leq C.C_{\beta }
\end{equation*}%
which completes the proof.

\textbf{Remark 3.1.} If $E_{1}$ does not have property ($\alpha $), we can
use another decomposition of $R^{n}$ to get an unconditional Schauder
decomposition of $L_{p,\gamma }\left( R^{n};E_{1}\right) $. But without
property ($\alpha $) we have to impose stronger conditions on $M$ to get $%
L_{p,\gamma }$ boundedness of the corresponding multiplier operator.

\textbf{Proof of Theorem 3.2. }For $j\in \mathbb{Z}$, let $s\left( j\right)
\in \mathbb{Z}$ and $t=t\left( j\right) \in \left\{ 1,2,...n\right\} $ be
the unique numbers satisfying $j=ns+t$. Set%
\begin{equation*}
D_{j}=\left( 0,2^{s}\right) ^{t-1}\times \left[ 2^{s},\right. \left.
2^{s+1}\right) \times \left( 0,2^{s}\right) ^{n-s}
\end{equation*}%
and define $\Delta _{j}=F^{-1}\chi _{D_{j}}F.$ Let $j=ns+t$\ be the unique
representation of $j.$ For $k\in \mathbb{Z}$, $r\in \mathbb{Z}^{n}$ with $%
0\leq r_{i}\leq 2^{k}-1$ define the operator $M_{j,r}^{k}$ by%
\begin{equation*}
M_{j,r}^{k}=M\left( y_{1},y_{2},...,y_{n}\right) ,
\end{equation*}%
where 
\begin{equation*}
\begin{array}{c}
y_{i}=r_{i}2^{s+1-k}\text{ for }i\in \left\{ 1,2,...,t-1\right\} \\ 
y_{t}=2^{s}+r_{t}2^{s-k}, \\ 
y_{i}=r_{i}2^{s-k},\text{ }i\in \left\{ t+1,t+2,...,n\right\} .%
\end{array}%
\end{equation*}

Then, by reasoning as the proof of Theorem 3.7 in $\left[ 9\right] $ we get
the assertion.

\begin{center}
\bigskip \textbf{4. Embeding theorems in Sobolev-Lions type spaces}
\end{center}

The embedding of Sobolev-Lions spaces play important roll in the regularity
theory of PDE with operator coefficients. In this section, we show
continuity of embedding operators in anisotropic Sobolev-Lions spaces.

\bigskip Let 
\begin{equation*}
X=L_{p,\gamma }\left( R^{n};E\right) \text{, }Y=W_{p,\gamma }^{l}\left(
R^{n};E\left( A\right) ,E\right) ,\ 
\end{equation*}%
\begin{equation*}
\text{ }l=\left( l_{1},l_{2},...,l_{n}\right) ,\text{ }\alpha =\left( \alpha
_{1},\alpha _{2},...,\alpha _{n}\right) ,\text{ }\varkappa =\left\vert
\alpha :l\right\vert =\dsum\limits_{k=1}^{n}\frac{\alpha _{k}}{l_{k}},\text{ 
}
\end{equation*}%
\begin{equation*}
\xi =\left( \xi _{1},\xi _{2},...,\xi _{n}\right) \in R^{n},\text{ }%
\left\vert \xi \right\vert ^{\alpha }=\prod\limits_{k=1}^{n}\left\vert \xi
_{k}\right\vert ^{\alpha _{k}}.
\end{equation*}%
From $\left[ \text{22, Lemma 3.1}\right] $ we have

\textbf{Lemma 4.1}. Assume $A$ is a $\varphi -$ positive linear operator on
a Banach space $E$. Then for any\ $h>0$ and $0\leq \mu \leq 1-\varkappa $
the operator-function \ 

\begin{equation*}
\Psi \left( \xi \right) =\Psi _{h}\left( \xi \right) =\left\vert \xi
\right\vert ^{\alpha }A^{1-\varkappa -\mu }h^{-\mu }\left[
A+\sum\limits_{k=1}^{n}\left\vert \xi _{k}\right\vert ^{l_{k}}+h^{-1}\right]
^{-1}
\end{equation*}%
is bounded in\ $E$ uniformly with respect to $\xi \in R^{n}$ and$\ h>0$ i.e.
there exists a constant $C_{\mu }$ such that 
\begin{equation}
\ \ \ \left\Vert \Psi _{h}\left( \xi \right) \right\Vert _{B\left( E\right)
}\leq C_{\mu }  \tag{4.1}
\end{equation}%
for all $\xi \in R^{n}$ and $h>0.$

One of main result of this section is the following:

\ \ \ \ \ \ \textbf{Theorem 4.1}. Let $\gamma \in A_{p}$ for $p\in \left(
1,\infty \right) .$ Assume $E$ is an UMD\ space and $A$ is a $\varphi -$
positive operator in $E$. Then for $0\leq \mu \leq 1-\varkappa $\ the
embedding 
\begin{equation*}
D^{\alpha }Y\subset L_{p\mathbf{,\gamma }}\left( R^{n};E\left(
A^{1-\varkappa -\mu }\right) \right)
\end{equation*}%
is a continuous and there exists a constant $C_{\mu }$ \ $>0$ depending only
on $\mu ,$ $p$, $\gamma $ such that 
\begin{equation}
\left\Vert D^{\alpha }u\right\Vert _{L_{p\mathbf{,\gamma }}\left(
R^{n};E\left( A^{1-\varkappa -\mu }\right) \right) }\leq C_{\mu }\left[
h^{\mu }\left\Vert u\right\Vert _{Y}+h^{-\left( 1-\mu \right) }\left\Vert
u\right\Vert _{X}\right]  \tag{4.2}
\end{equation}%
for $u\in Y$ and $0<h\leq h_{0}<\infty .$

\ \ \textbf{Proof}. It is clear to see that

\begin{equation}
A^{1-\alpha -\mu }D^{\alpha }u=F^{-1}FA^{1-\varkappa -\mu }D^{\alpha
}u=F^{-1}A^{1-\varkappa -\mu }FD^{\alpha }u=  \tag{4.3}
\end{equation}%
\begin{equation*}
F^{-1}A^{1-\varkappa -\mu }\left( i\xi \right) ^{\alpha }Fu=F^{-1}\left(
i\xi \right) ^{\alpha }A^{1-\varkappa -\mu }Fu.
\end{equation*}%
Hence, denoting \ $Fu$ by $\hat{u},$ we get from $\left( 4.3\right) $ the
following estimate%
\begin{equation*}
C_{2}\left\Vert F^{-1}\left( i\xi \right) ^{\alpha }A^{1-\varkappa -\mu }%
\hat{u}\right\Vert _{X}\leq \left\Vert D^{\alpha }u\right\Vert _{L_{p,\gamma
}\left( R^{n};E\left( A^{1-\varkappa -\mu }\right) \right) }\leq
\end{equation*}%
\begin{equation*}
C_{1}\left\Vert F^{-1}\left( i\xi \right) ^{\alpha }A^{1-\varkappa -\mu }%
\hat{u}\right\Vert _{X},
\end{equation*}%
where $C_{1}$, $C_{2}$ are positive constants depending only of $p$ and $%
\gamma $. Similarly, there exist positive constants $M_{1}$ and $M_{2}$ such
that for $u\in Y$ we have 
\begin{equation*}
M_{1}\left\Vert u\right\Vert _{Y}\leq \left\Vert F^{-1}\hat{u}\right\Vert
_{X}+\sum\limits_{k=1}^{n}\left\Vert F^{-1}\left[ \left( i\xi _{k}\right)
^{l_{k}}\hat{u}\right] \right\Vert _{X}\leq M_{2}\left\Vert u\right\Vert
_{Y}.
\end{equation*}%
Therefore, for proving the inequality $\left( 4.2\right) $ it suffices to
show 
\begin{equation*}
\left\Vert F^{-1}\left( i\xi \right) ^{\alpha }A^{1-\varkappa -\mu }\hat{u}%
\right\Vert _{X}\leq
\end{equation*}%
\begin{equation}
C_{\mu }(h^{\mu }\left\Vert F^{-1}A\hat{u}\right\Vert
_{X}+\sum\limits_{k=1}^{n}\left\Vert F^{-1}\left[ \left( i\xi _{k}\right)
^{l_{k}}\hat{u}\right] \right\Vert _{X}+h^{-\left( 1-\mu \right) }\left\Vert
F^{-1}\hat{u}\right\Vert _{X}).  \tag{4.4}
\end{equation}%
Therefore, the inequality $\left( 4.2\right) $ will follow if we prove the
following estimate \ 
\begin{equation}
\left\Vert F^{-1}\left[ \xi ^{\alpha }A^{1-\varkappa -\mu }\hat{u}\right]
\right\Vert _{X}\leq C_{\mu }\left\Vert F^{-1}G\left( \xi \right) \hat{u}%
\right\Vert _{X}\text{.}  \tag{4.5}
\end{equation}%
for $u\in Y,$ where%
\begin{equation*}
G\left( \xi \right) =h^{\mu }\left[ A+\sum\limits_{k=1}^{n}\left\vert \xi
_{k}\right\vert ^{l_{k}}+h^{-\left( 1-\mu \right) }\right] .
\end{equation*}%
\ Due to positivity of\ $A,$ the operator function $G\left( \xi \right) $
has a bounded inverse in\ $E$ for all $\xi \in R^{n}$ and $h>0.$ So, we can
set%
\begin{equation}
F^{-1}\xi ^{\alpha }A^{1-\varkappa -\mu }\hat{u}=F^{-1}\xi ^{\alpha
}A^{1-\varkappa -\mu }G^{-1}\left( \xi \right) \left[ h^{\mu }\left(
A+\sum\limits_{k=1}^{n}\left\vert \xi _{k}\right\vert ^{l_{k}}\right)
+h^{-\left( 1-\mu \right) }\right] \hat{u}.  \tag{4.6}
\end{equation}%
The inequality $\left( 4.5\right) $ will follow immediately from $\left(
4.6\right) $ if we can prove that the operator-function $\Psi _{h}=\xi
^{\alpha }A^{1-\varkappa -\mu }G^{-1}\left( \xi \right) $ is a multiplier in 
$M_{p,\gamma }^{q,\gamma }\left( E\right) $ uniformly with respect to\ $h.$
So, by Theorem 3.1 it suffices to show that the set%
\begin{equation*}
B\left( \xi ,h\right) =\left\{ \xi ^{\beta }D^{\beta }\Psi _{h}\left( \xi
\right) ;\text{ }\xi \in R^{n}\backslash \left\{ 0\right\} ,\text{ }\beta
_{j}\in \left\{ 0,1\right\} \right\}
\end{equation*}%
\ is $R-$bounded uniformly in $h$, i.e.%
\begin{equation}
\sup\limits_{h}R\left\{ B\left( \xi ,h\right) \right\} \leq M.  \tag{4.7}
\end{equation}%
By Lemma 4.1 there exists a constant \ $C_{\mu }>0$ such that the following
uniform estimate holds%
\begin{equation}
\left\Vert \Psi _{h}\left( \xi \right) \right\Vert _{B\left( E\right) }\leq
C_{\mu }.  \tag{4.8}
\end{equation}%
\ Let first, $\beta =\left( \beta _{1},...\beta _{n}\right) $ where \ $\beta
_{k}=1$ and $\beta =0$ for $j\neq k$. Then, by using the resolvent
properties of $A$ we obtain

\begin{equation*}
\left\vert \frac{\partial }{\partial \xi _{k}}\Psi _{h}\left( \xi \right)
\right\vert \leq \prod\limits_{k=1}^{n}\left( i\right) ^{\left\vert \alpha
\right\vert }\alpha _{k}\left\vert \xi _{1}^{\alpha _{1}}...\xi
_{k-1}^{\alpha _{k-1}}\xi _{k}^{\alpha _{k}-1}...\xi _{n}^{\alpha
_{n}}\right\vert
\end{equation*}%
\begin{equation*}
\left\Vert A^{1-\varkappa -\mu }\left[ h^{\mu }\left(
A+\sum\limits_{k=1}^{n}\left\vert \xi _{k}\right\vert ^{l_{k}}\right)
+h^{-\left( 1-\mu \right) }\right] ^{-1}\right\Vert +
\end{equation*}%
\begin{equation*}
\left\vert \xi \right\vert ^{\alpha }\left\Vert A^{1-\varkappa -\mu }\left[
h^{\mu }\left( A+\sum\limits_{k=1}^{n}\left\vert \xi _{k}\right\vert
^{l_{k}}\right) +h^{-\left( 1-\mu \right) }\right] ^{-2}\right\Vert
h\left\vert \xi _{k}\right\vert ^{l_{k}-1}\leq
\end{equation*}%
\begin{equation*}
C_{\mu }\left\vert \xi _{k}\right\vert ^{-1},\text{ }k=1,2...n.
\end{equation*}%
Repeating the above process, we obtain that there exists a constant \ $%
C_{\mu }>0$ \ depending only $\mu $ such that%
\begin{equation*}
\ \left\vert \xi ^{\beta }\right\vert \ \left\Vert D^{\beta }\Psi _{h}\left(
\xi \right) \right\Vert _{B\left( E\right) }\leq C_{\mu }
\end{equation*}%
for $\beta =\left( \beta _{1},...\beta _{n}\right) $, $\beta _{k}\in \left\{
0,1\right\} $ and for all $\xi \in R^{n}$, $\xi \neq 0$. \ Due to $R$%
-positivity of $A$ and by $\left( 4.9\right) $ we obtain that the set 
\begin{equation*}
B_{0}\left( \xi \right) =\left\{ AG^{-1}\left( \xi ,h\right) ;\text{ }\xi
\in R^{n}\backslash \left\{ 0\right\} \text{, }\beta _{j}\in \left\{
0,1\right\} \right\}
\end{equation*}%
is $R$ bounded uniformly in $h$. Then, by virtue of Kahane's contraction
principle $\left[ \text{8, Lemma 3.5}\right] $ and by $\left( 4.9\right) $\
we obtain that the set%
\begin{equation*}
B_{1}\left( \xi ,h\right) =\left\{ AD^{-2}\left( \xi ,h\right) ;\text{ }\xi
\in R^{n}\backslash \left\{ 0\right\} \text{, }\beta _{j}\in \left\{
0,1\right\} \right\}
\end{equation*}%
is uniformly $R$-bounded. Moreover, by using the inequalities of moment for
positive operators and Young's we get that 
\begin{equation}
\left\Vert \Psi _{h}\left( \xi \right) u\right\Vert \leq C_{\mu }\left(
\left\Vert Au\right\Vert +\dsum\limits_{k=1}^{n}\left\vert \xi
_{k}\right\vert ^{l_{k}}\left\Vert u\right\Vert \right) ,  \tag{4.10}
\end{equation}%
where 
\begin{equation*}
u=G^{-1}\left( \xi ,h\right) f\text{, }f\in E.
\end{equation*}

Then thanks to $R$-boundedness of $B_{i}\left( \xi ,\lambda \right) $\ we
have 
\begin{equation}
\int\limits_{0}^{1}\left\Vert \sum\limits_{j=1}^{m}r_{j}\left( y\right)
B_{i}\left( \eta _{j},h\right) u_{j}\right\Vert _{E}dy\leq
C\int\limits_{0}^{1}\left\Vert \sum\limits_{j=1}^{m}r_{j}\left( y\right)
u_{j}\right\Vert _{E}dy,  \tag{4.11}
\end{equation}%
for all $\xi _{1},\xi _{2},...,\xi _{m}\in R^{n}$, $\eta _{j}=\left( \xi
_{j1,}\xi _{j2},...,\xi _{jn}\right) \in R^{n},$ $u_{1,}u_{2},...,u_{m}\in E$%
, $m\in \mathbb{N}$, where $\left\{ r_{j}\right\} $ is a sequence of
independent symmetric $\left\{ -1,1\right\} $-valued random variables on $%
\left[ 0,1\right] $. Thus, in view of Kahane's contraction principle,
additional and product properties of $R$-bounded operators and $\left(
4.10\right) $, $\left( 4.11\right) $ we obtain 
\begin{equation}
\int\limits_{0}^{1}\left\Vert \sum\limits_{j=1}^{m}r_{j}\left( y\right) \Psi
\left( \eta _{j},h\right) u_{j}\right\Vert _{E}dy\leq
C\int\limits_{0}^{1}\sum\limits_{i=0}^{1}\left\Vert
\sum\limits_{j=1}^{m}B_{i}\left( \eta _{j},h\right) r_{j}\left( y\right)
u_{j}\right\Vert _{E}dy\leq  \tag{4.12}
\end{equation}

\begin{equation*}
C\int\limits_{0}^{1}\left\Vert \sum\limits_{j=1}^{m}r_{j}\left( y\right)
u_{j}\right\Vert _{E}dy.
\end{equation*}%
The estimate $\left( 4.12\right) $ implies $R$-boundedness of the set $%
B\left( \xi ,h\right) $, which implies the assertion$.$

\ \ It is possible to state Theorem 4.1 in a more general setting. For this
aim, we use the concept of extension operator.

\ \textbf{Condition 4.1}. Let $\gamma \in A_{p}$ for $p\in \left( 1,\infty
\right) $. Let $A$ be a positive operator in UMD space $E.$\ Assume a region
\ $\Omega \subset R^{n}$ such that there exists bounded linear extension
operator$\ B$ from $W_{p,\gamma }^{l}\left( \Omega ,E\left( A\right)
,E\right) $ to $Y$ for $1<p<\infty .$

\textbf{Remark 4.1}. If \ $\Omega \subset R^{n}$ is a region satisfying the
strong \ $l-$horn condition \ (see $\left[ 3\right] $, p.117 for $E=\mathbb{C%
},$ $A$ $=I$ \ and $\gamma \left( x\right) \equiv 1$)\ then for $1<p<\infty $
there exists a bounded linear extension operator from \ $W_{p}^{l}\left(
\Omega \right) =W_{p}^{l}\left( \Omega ;\mathbb{C},\mathbb{C}\right) $ to $%
W_{p}^{l}\left( R^{n}\right) =W_{p}^{l}\left( R^{n};\mathbb{C},\mathbb{C}%
\right) .$

\textbf{Theorem 4.2}. Assume conditions of Theorem 4.1 and Condition 4.1 are
satisfied. Then for $0\leq \mu \leq 1-\varkappa $ the embedding 
\begin{equation*}
D^{\alpha }W_{p,\gamma }^{l}\left( \Omega ;E\left( A\right) ,E\right)
\subset L_{p,\gamma }\left( \Omega ;E\left( A^{1-\varkappa -\mu }\right)
\right)
\end{equation*}%
is continuous and there exists a constant $C_{\mu }$ depending only of $\mu
, $ $p$, $\gamma $ such that 
\begin{equation}
\left\Vert D^{\alpha }u\right\Vert _{L_{p,\gamma }\left( \Omega ;E\left(
A^{1-\varkappa -\mu }\right) \right) }\leq  \tag{4.10}
\end{equation}%
\begin{equation*}
C_{\mu }\left[ h^{\mu }\left\Vert u\right\Vert _{W_{p,\gamma }^{l}\left(
\Omega ;E\left( A\right) ,E\right) }+h^{-\left( 1-\mu \right) }\left\Vert
u\right\Vert _{L_{p,\gamma }\left( \Omega ;E\right) }\right]
\end{equation*}%
for $u\in W_{p,\gamma }^{l}\left( \Omega ;E\left( A\right) ,E\right) $ and $%
0<h\leq h_{0}<\infty .$

\ \textbf{Proof}.\ It is suffices to prove the estimate $\left( 4.10\right)
. $ Let $B$ is a bounded linear extension operator\ from \ $W_{p,\gamma
}^{l}\left( \Omega ;E\left( A\right) ,E\right) $\ to $W_{p,\gamma
}^{l}\left( R^{n};E\left( A\right) ,E\right) $ and let $B_{\Omega }$\ be\
the restriction operator from \ $R^{n}$ to $\Omega .$ Then for any $u\in
W_{p,\gamma }^{l}\left( \Omega ;E\left( A\right) ,E\right) $\ we have 
\begin{equation*}
\left\Vert D^{\alpha }u\right\Vert _{L_{p,\gamma }\left( \Omega ;E\left(
A^{1-\varkappa -\mu }\right) \right) }=\left\Vert D^{\alpha }B_{\Omega
}Bu\right\Vert _{L_{p,\gamma }\left( \Omega ;E\left( A^{1-\varkappa -\mu
}\right) \right) }\leq
\end{equation*}%
\begin{equation*}
C_{\mu }\left[ h^{\mu }\left\Vert Bu\right\Vert _{W_{p,\gamma }^{l}\left(
R^{n};E\left( A\right) ,E\right) }+h^{-\left( 1-\mu \right) }\left\Vert
Bu\right\Vert _{L_{p,\gamma }\left( R^{n};E\right) }\right] \leq
\end{equation*}%
\begin{equation*}
C_{\mu }\left[ h^{\mu }\left\Vert u\right\Vert _{W_{p,\gamma }^{l}\left(
\Omega ;E\left( A\right) E\right) }+h^{-\left( 1-\mu \right) }\left\Vert
u\right\Vert _{L_{p,\gamma }\left( \Omega ;E\right) }\right] .
\end{equation*}%
\textbf{Result 4.1}. Assume the conditions of Theorem 4.2 are satisfied.
Then for $u\in W_{p,\gamma }^{l}\left( \Omega ;E\left( A\right) ,E\right) $
we have the following multiplicative estimate 
\begin{equation}
\left\Vert D^{\alpha }u\right\Vert _{L_{p,\gamma }\left( \Omega ;E\left(
A^{1-\varkappa -\mu }\right) \right) }\leq C_{\mu }\left\Vert u\right\Vert
_{W_{p,\gamma }^{l}\left( \Omega ;E\left( A\right) ,E\right) }^{1-\mu
}.\left\Vert u\right\Vert _{L_{p,\gamma }\left( \Omega ;E\right) }^{\mu }. 
\tag{4.11}
\end{equation}%
\ Indeed, setting 
\begin{equation*}
h=\left\Vert u\right\Vert _{L_{p,\gamma }\left( \Omega ;E\right)
}.\left\Vert u\right\Vert _{W_{p,\gamma }^{l}\left( \Omega ;E\left( A\right)
,E\right) }^{-1}
\end{equation*}%
\ in\ $\left( 4.10\right) $ we obtain $\left( 4.11\right) .$

\bigskip \textbf{Theorem 4}.\textbf{3.} Suppose conditions of Theorem 4.1
are hold. Then for $0<\mu <1-\varkappa $ the embedding 
\begin{equation*}
D^{\alpha }Y\subset L_{p,\gamma }\left( R^{n};\left( E\left( A\right)
,E\right) _{\varkappa +\mu ,p}\right)
\end{equation*}%
is continuous and there exists a constant $C_{\mu }$ depending only of $\mu
, $ $p$, $\gamma $ such that

\begin{equation}
\left\Vert D^{\alpha }u\right\Vert _{L_{p,\gamma }\left( R^{n};\left(
E\left( A\right) ,E\right) _{\varkappa +\mu ,p}\right) }\leq h^{\mu
}\left\Vert u\right\Vert _{Y}+h^{-\left( 1-\mu \right) }\left\Vert
u\right\Vert _{X}  \tag{4.12}
\end{equation}%
for $u\in Y$ and $0<h\leq h_{0}<\infty .$

\textbf{Proof. }It is sufficient to prove the estimate $\left( 4.12\right) $
for $u\in Y.$ By definition of interpolation spaces $\left( E\left( A\right)
,E\right) _{\varkappa +\mu ,p}$ (see $\left[ \text{27, \S 1.14.5}\right] $)
the estimate $\left( 4.12\right) $ is equivalent to the inequality

\begin{equation}
\left\Vert F^{-1}y^{1-\varkappa -\mu -\frac{1}{p}}\left[ A^{\chi +\mu
}\left( A+y\right) ^{-1}\right] \xi ^{\alpha }\hat{u}\right\Vert
_{L_{p,\gamma }\left( R_{+}^{n+1};E\right) }  \tag{4.13}
\end{equation}

\begin{equation*}
\leq C_{\mu }\left\Vert F^{-1}\left[ h^{\mu
}(A+\sum\limits_{k=1}^{n}A+\sum\limits_{k=1}^{n}\left\vert \xi
_{k}\right\vert ^{l_{k}}+h^{-\left( 1-\mu \right) }\right] \hat{u}%
\right\Vert _{L_{p,\gamma }\left( R^{n};E\right) }.
\end{equation*}%
By multiplier properties, the inequality $\left( 4.13\right) $ will follow
immediately if we will prove that the operator-function

\begin{equation*}
\Psi =\left( i\xi \right) ^{\alpha }y^{1-\varkappa -\mu -\frac{1}{p}}A^{\chi
+\mu }\left( A+y\right) ^{-1}\left[ h^{\mu }\left(
A+\sum\limits_{k=1}^{n}\left\vert \xi _{k}\right\vert ^{l_{k}}\right)
+h^{-\left( 1-\mu \right) }\right] ^{-1}
\end{equation*}%
is a multiplier from $X$ to $L_{p,\gamma }\left( R^{n};L_{p}\left(
R_{+};E\right) \right) .$ This fact is proved by the same manner as Theorem
4.1. Therefore, we get the estimate $\left( 4.12\right) .$

In a similar way, as the Theorem 4.2 we obtain

\textbf{Theorem 4.4}. Suppose conditions of Theorem 4.2 are hold$.$ Then for 
$0<\mu <1-\varkappa $ the embedding 
\begin{equation*}
D^{\alpha }W_{p,\gamma }^{l}\left( \Omega ;E\left( A\right) ,E\right)
\subset L_{p,\gamma }\left( \Omega ;\left( E\left( A\right) ,E\right)
_{\varkappa +\mu ,p}\right)
\end{equation*}%
is continuous and there exists a constant $C_{\mu }$ depending only of $\mu
, $ $p$, $\gamma $ such that 
\begin{equation}
\left\Vert D^{\alpha }u\right\Vert _{L_{p,\gamma }\left( \Omega ,\left(
E\left( A\right) ,E\right) _{\varkappa +\mu ,p}\right) }\leq C_{\mu }\left[
h^{\mu }\left\Vert u\right\Vert _{W_{p,\gamma }^{l}\left( \Omega ;E\left(
A\right) ,E\right) }+h^{-\left( 1-\mu \right) }\left\Vert u\right\Vert
_{L_{p,\gamma }\left( \Omega ;E\right) }\right]  \tag{4.14}
\end{equation}%
for $u\in W_{p,\gamma }^{l}\left( \Omega ;E\left( A\right) ,E\right) $ and $%
0<h\leq h_{0}<\infty .$

\textbf{Result 4. 2}. Suppose the conditions of Theorem 4.2 are hold. Then
for $u\in W_{p,\gamma }^{l}\left( \Omega ;E\left( A\right) ,E\right) $\ we
have the following multiplicative estimate 
\begin{equation}
\left\Vert D^{\alpha }u\right\Vert _{L_{p,\gamma }\left( \Omega ;\left(
E\left( A\right) ,E\right) _{\varkappa +\mu ,p}\right) }\leq C_{\mu
}\left\Vert u\right\Vert _{W_{p,\gamma }^{l}\left( \Omega ;E\left( A\right)
,E\right) }^{1-\mu }\left\Vert u\right\Vert _{L_{p,\gamma }\left( \Omega
;E\right) }^{\mu }.  \tag{4.15}
\end{equation}%
Indeed setting $h=\left\Vert u\right\Vert _{L_{p,\gamma }\left( \Omega
;E\right) }.\left\Vert u\right\Vert _{W_{p,\gamma }^{l}\left( \Omega
;E\left( A\right) ,E\right) }^{-1}$ in $\left( 4.14\right) $ we obtain\ $%
\left( 4.15\right) .$

From Theorem 4.2 we obtain

\textbf{Result 4.3}. Assume the conditions of Theorem 4.2 are satisfied for $%
l_{1}=l_{2}=\ldots =l_{n}=m.$ Then for $0\leq \mu \leq 1-\varkappa $ the
embedding 
\begin{equation*}
D^{\alpha }W_{p,\gamma }^{m}\left( \Omega ;E\left( A\right) ,E\right)
\subset L_{p,\gamma }\left( \Omega ;E\left( A^{1-\varkappa -\mu }\right)
\right)
\end{equation*}%
is continuous and there exists a constant $C_{\mu }$ depending only of $\mu
, $ $p$, $\gamma $ such that 
\begin{equation*}
\left\Vert D^{\alpha }u\right\Vert _{L_{p,\gamma }\left( \Omega ;E\left(
A^{1-\varkappa -\mu }\right) \right) }\leq C_{\mu }\left[ h^{\mu }\left\Vert
u\right\Vert _{W_{p,\gamma }^{m}\left( \Omega ;E\left( A\right) ,E\right)
}+h^{-\left( 1-\mu \right) }\left\Vert u\right\Vert _{L_{p,\gamma }\left(
\Omega ;E\right) }\right]
\end{equation*}%
for $u\in W_{p,\gamma }^{m}\left( \Omega ;E\left( A\right) ,E\right) $ and\ $%
0<h\leq h_{0}<\infty $ where 
\begin{equation*}
\varkappa =\frac{\left\vert \alpha \right\vert }{m}.
\end{equation*}

\textbf{Result 4.3}. If $E=H$, where $H$ is a Hilbert space and \ $%
p_{k}=q_{k}=2,$ $\Omega =\left( 0,T\right) ,$ $l_{1}=l_{2}=\ldots =l_{n}=m$ $%
,$ $A=A^{\times }\geq c^{2}I,$\ $\gamma \left( x\right) \equiv 1$\ then we
obtain the well known Lions-Peetre $\left[ \text{14}\right] $ result. \
Moreover, the result of Lions-Peetre is improving even in the one
dimensional case and for non selfedjoint positive operators $A.$

From Theorems 4.2 we obtain

\textbf{Result 4.4}. Suppose the conditions of Theorem 4.2 are satisfied for 
$\gamma \left( x\right) \equiv 1.$ Then for $0\leq \mu \leq 1-\varkappa $
the embedding 
\begin{equation*}
D^{\alpha }W_{p}^{l}\left( \Omega ;E\left( A\right) ,E\right) \subset
L_{p}\left( \Omega ;E\left( A^{1-\varkappa -\mu }\right) \right)
\end{equation*}%
is continuous and there exists a constant $C_{\mu }$ \ $>0$ depending only
of $\mu ,$ $p$, $\gamma $ such that 
\begin{equation*}
\left\Vert D^{\alpha }u\right\Vert _{L_{p}\left( \Omega ;E\left(
A^{1-\varkappa -\mu }\right) \right) }\leq C_{\mu }\left[ h^{\mu }\left\Vert
u\right\Vert _{W_{p}^{l}\left( \Omega ;E\left( A\right) ,E\right)
}+h^{-\left( 1-\mu \right) }\left\Vert u\right\Vert _{L_{p}\left( \Omega
;E\right) }\right]
\end{equation*}%
for $u\in W_{p}^{l}\left( \Omega ;E\left( A\right) ,E\right) $ and $0<h\leq
h_{0}<\infty $.

Moreover, if $\Omega $ is a bounded domain in $R^{n}$ and $A^{-1}$ is a
compact operator in $E,$ then for $0<\mu \leq 1-\varkappa $ the embedding%
\begin{equation*}
D^{\alpha }W_{p}^{l}\left( \Omega ;E\left( A\right) ,E\right) \subset
L_{p}\left( \Omega ;E\left( A^{1-\varkappa -\mu }\right) \right)
\end{equation*}%
is compact.

If $E=\mathbb{C},$ $A=I$, $\gamma \left( x\right) \equiv 1$ we get the
embedding $D^{\alpha }W_{p}^{l}\left( \Omega \right) \subset L_{p}\left(
\Omega \right) $ proved in $\left[ 3\right] $ for Sobolev spaces $%
W_{p}^{l}\left( \Omega \right) .$

\bigskip Let\ $s>0.$ Consider the following sequence space (see e.g. $\left[ 
\text{27, \S\ 1.18}\right] $)%
\begin{equation*}
l_{q}^{s}=\left\{ u=\left\{ u_{i}\right\} ,\text{ }i=1,2,...,\infty ,\text{ }%
u_{i}\in \mathbb{C}\right\}
\end{equation*}%
\ with the norm$\ $%
\begin{equation*}
\ \ \ \ \ \ \left\Vert u\right\Vert _{l_{q}^{s}}=\left(
\sum\limits_{i=1}^{\infty }2^{i\nu s}\left\vert u_{i}\right\vert ^{p}\right)
^{\frac{1}{q}}<\infty ,\text{ }\nu \in \left( 1,\infty \right) .
\end{equation*}%
Note that, \ $l_{q}^{0}=l_{q}.$ Let $A$ be infinite matrix defined in $%
l_{\nu }$ such that \ $D\left( A\right) =l_{q}^{s},$ $A=\left[ \delta
_{ij}2^{si}\right] ,\ $where \ $\delta _{ij}=0$, when \ $i\neq j,$ \ $\delta
_{ij}=1,$ when $i=j=1,2,...,\infty .$

It is clear to see that the operator $A$\ is positive in\ $l_{q}$. From
Theorem 4.2 we obtain the following results:

\textbf{Result 4.5. }Suppose the conditions of Theorem 4.2 are satisfied for 
$E=\mathbb{C}$. Then for $0\leq \mu \leq 1-\varkappa $, $1<p<\infty $ the
embedding 
\begin{equation*}
D^{\alpha }W_{p,\gamma }^{l}\left( \Omega ,l_{q}^{s},l_{q}\right) \subset
L_{p,\gamma }\left( \Omega ,l_{q}^{s\left( 1-\varkappa -\mu \right) }\right)
\end{equation*}%
is continuous and there exists a constant $C_{\mu }$ \ $>0$ depending only
of $\mu ,$ $p$, $q,$ $\gamma $ such that 
\begin{equation*}
\left\Vert D^{\alpha }u\right\Vert _{L_{p\mathbf{,\gamma }}\left( \Omega
;l_{q}^{s\left( 1-\varkappa -\mu \right) }\right) }\leq C_{\mu }\left[
h^{\mu }\left\Vert u\right\Vert _{W_{p,\gamma }^{l}\left( \Omega ;l_{\nu
}^{s},l_{\nu }\right) }+h^{-\left( 1-\mu \right) }\left\Vert u\right\Vert
_{L_{p\mathbf{,\gamma }}\left( \Omega ;l_{\nu }\right) }\right]
\end{equation*}%
for\ $u\in W_{p,\gamma }^{l}\left( \Omega ,l_{q}^{s},l_{q}\right) $ and $%
0<h\leq h_{0}<\infty .$

\textbf{Result 4.6. }Suppose the conditions of Theorem 4.2 are hold for $E=%
\mathbb{C}$. Then for $0<\mu \leq 1-\varkappa ,$ $1<p<\infty $ the embedding 
\begin{equation*}
D^{\alpha }W_{p,\gamma }^{l}\left( \Omega ,l_{q}^{s},l_{q}\right) \subset
L_{p,\gamma }\left( \Omega ,l_{q}^{s\left( 1-\varkappa -\mu \right) }\right)
\end{equation*}%
is compact.

\textbf{Result 4.7. }For $0\leq \mu \leq 1-\varkappa $, $1<p<\infty $ the
embedding 
\begin{equation*}
D^{\alpha }W_{p}^{l}\left( \Omega ,l_{q}^{s},l_{q}\right) \subset
L_{p}\left( \Omega ,l_{q}^{s\left( 1-\varkappa -\mu \right) }\right)
\end{equation*}%
is a continuous and there exists a constant $C_{\mu }$ \ $>0$, depending
only of $\mu ,$ $p$, $q,$ $\gamma $ such that 
\begin{equation*}
\left\Vert D^{\alpha }u\right\Vert _{L_{p}\left( \Omega ;l_{q}^{s\left(
1-\varkappa -\mu \right) }\right) }\leq C_{\mu }\left[ h^{\mu }\left\Vert
u\right\Vert _{W_{p}^{l}\left( \Omega ;l_{\nu }^{s},l_{\nu }\right)
}+h^{-\left( 1-\mu \right) }\left\Vert u\right\Vert _{L_{p}\left( \Omega
;l_{q}\right) }\right]
\end{equation*}%
for\ $u\in W_{p}^{l}\left( \Omega ,l_{q}^{s},l_{q}\right) $ and $0<h\leq
h_{0}<\infty .$

Note that, these results haven't been obtained with classical method until
now.

\bigskip \textbf{\ }

\begin{center}
\bigskip \textbf{5. Separable differential operators in weighted Lebesque
spaces}
\end{center}

Firstly, consider the leading part of the equation $\left( 1.1\right) $,
i.e. consider the following equation 
\begin{equation}
\ L_{0}u=\sum\limits_{\left\vert \alpha \right\vert =2l}a_{\alpha }D^{\alpha
}u+Au+\lambda u=f,  \tag{5.1}
\end{equation}%
where $a_{\alpha }$ are complex numbers, $l\in \mathbb{N},$ $A$ is a linear
operator in a Banach space $E$ and $\lambda $ is a complex parameter.

Let 
\begin{equation*}
X=L_{p,\gamma }\left( R^{n};E\right) ,\text{ }Y=W_{p,\gamma }^{2l}\left(
R^{n};E\left( A\right) ,E\right) .
\end{equation*}

\textbf{Condition 5.1. }\ Let 
\begin{equation*}
\text{(a) }K\left( \xi \right) =\sum\limits_{\left\vert \alpha \right\vert
=2l}a_{\alpha }\left( i\xi _{1}\right) ^{\alpha _{1}}\left( i\xi _{2}\right)
^{\alpha _{2}}...\left( i\xi _{n}\right) ^{\alpha _{n}}\in S\left( \varphi
_{1}\right)
\end{equation*}%
for $0\leq \varphi _{1}<\pi ;$

(b) There exists the positive constat $M_{0}$ so that%
\begin{equation*}
\text{ }\left\vert K\left( \xi \right) \right\vert \geq
M_{0}\dsum\limits_{k=1}^{n}\xi _{k}^{2l}\text{ for all }\xi \in R^{n},\text{ 
}\xi \neq 0.
\end{equation*}

In this section we prove the following result

\textbf{Theorem 5.1.}\ Suppose the following conditions hold:

(1) Condition 5.1 is hold;

(2) $\gamma \in A_{p}$ for $p\in \left[ 1,\infty \right] $;

(3) $A$ is a $R-$positive operator in UMD\ space $E$ for $0\leq \varphi <\pi
-\varphi _{1}$.

Then for all $\ f\in X$ and $\lambda \in S\left( \varphi _{1}\right) $
equation $\left( 6.1\right) $ has an unique solution $u$ that belongs to
space $Y$ and the coercive uniform estimate holds

\begin{equation}
\sum\limits_{\left\vert \alpha \right\vert \leq 2l}\left\vert \lambda
\right\vert ^{1-\frac{\left\vert \alpha \right\vert }{2l}}\left\Vert
D^{\alpha }u\right\Vert _{X}+\left\Vert Au\right\Vert _{X}\leq C\left\Vert
f\right\Vert _{X}.  \tag{5.2}
\end{equation}

\ \textbf{Proof}. By applying the Fourier transform to the equation $\left(
5.1\right) $\ we get 
\begin{equation}
\left[ K\left( \xi \right) +A+\lambda \right] \hat{u}\left( \xi \right) =f^{%
\symbol{94}}\left( \xi \right) ,  \tag{5.3}
\end{equation}%
where 
\begin{equation*}
K\left( \xi \right) =\sum\limits_{\left\vert \alpha \right\vert
=2l}a_{\alpha }\left( i\xi _{1}\right) ^{\alpha _{1}}\left( i\xi _{2}\right)
^{\alpha _{2}}...\left( i\xi _{n}\right) ^{\alpha _{n}}.
\end{equation*}

Since\ $K\left( \xi \right) \in S\left( \varphi _{1}\right) $ for all\ $\xi
\in R^{n},$ the operator\ $A+K\left( \xi \right) $ is invertible in $E$. So,
we obtain that the solution of the equation $\left( 5.3\right) $ can be
represented in the form 
\begin{equation}
u\left( x\right) =F^{-1}\left[ A+K\left( \xi \right) +\lambda \right]
^{-1}f^{\symbol{94}}.  \tag{5.4}
\end{equation}%
By using $\left( 5.4\right) $ we have

\begin{equation*}
\left\Vert Au\right\Vert _{X}=\left\Vert F^{-1}A\left[ A+K\left( \xi \right)
+\lambda \right] ^{-1}f^{\symbol{94}}\right\Vert _{X},
\end{equation*}%
\begin{equation*}
\left\Vert D^{\alpha }u\right\Vert _{X}=\left\Vert F^{-1}\left( i\xi
_{1}\right) ^{\alpha _{1}}\left( i\xi _{2}\right) ^{\alpha _{2}}...\left(
i\xi _{n}\right) ^{\alpha _{n}}\left[ A+K\left( \xi \right) +\lambda \right]
^{-1}f^{\symbol{94}}\right\Vert _{X}.
\end{equation*}%
Hence, it is suffices to show that the operator-functions 
\begin{equation*}
\sigma _{1,\lambda }\left( \xi \right) =A\left[ A+K\left( \xi \right)
+\lambda \right] ^{-1}\text{, }
\end{equation*}%
\begin{equation*}
\sigma _{2,\lambda }\left( \xi \right) =\sum\limits_{\left\vert \alpha
\right\vert \leq 2l}\xi _{1}^{\alpha _{1}}\xi _{2}^{\alpha _{2}}...\xi
_{n}^{\alpha _{n}}\left\vert \lambda \right\vert ^{1-\frac{\left\vert \alpha
\right\vert }{2l}}\left[ A+K\left( \xi \right) +\lambda \right] ^{-1}
\end{equation*}%
are multipliers in $X.$ To see this, it is suffices to show that the
following collections 
\begin{equation*}
\left\{ \xi ^{\beta }D^{\beta }\sigma _{1,\lambda }\left( \xi \right) :\xi
\in R^{n}/\left\{ 0\right\} ,\text{ }\beta \in U_{n}\right\} ,\left\{ \xi
^{\beta }D^{\beta }\sigma _{2,\lambda }\left( \xi \right) :\xi \in
R^{n}/\left\{ 0\right\} ,\text{ }\beta \in U_{n}\right\}
\end{equation*}%
are $R-$bounded in $E$ uniformly in $\lambda $, where%
\begin{equation*}
U=\left\{ \beta =\beta _{1},...,\beta _{n}),\text{ }\beta _{i}\in \left\{
0,1\right\} \right\} \text{. }
\end{equation*}%
Due to $R-$positivity of $A$, the set 
\begin{equation*}
\left\{ \sigma _{1,\lambda }\left( \xi \right) :\xi \in R^{n}/\left\{
0\right\} ,\text{ }\beta \in U_{n}\right\}
\end{equation*}%
is $R$-bounded. Moreover, by using the same reasoning as used in the proof
of Theorem 4.1 and in view of (3) condition we obtain that the set%
\begin{equation*}
\left\{ \sigma _{2,\lambda }\left( \xi \right) :\xi \in R^{n}/\left\{
0\right\} ,\text{ }\beta \in U_{n}\right\}
\end{equation*}%
is $R$-bounded uniformly in $\lambda \in S\left( \varphi _{1}\right) $. Then
by virtue of of Kahane's contraction principle, by product properties of the
collection of $R$ -bounded operators (see e.g. Lemma 3.5., Proposition 3.4.
in $\left[ 8\right] $) and due to $R-$positivity of operator $A$ we obtain 
\begin{equation}
\sup\limits_{\lambda \in S\left( \varphi _{1}\right) }R\left\{ \xi ^{\beta
}D^{\beta }\sigma _{1,\lambda }\left( \xi \right) :\xi \in R^{n}/\left\{
0\right\} ,\text{ }\beta \in U_{n}\right\} \leq C,  \tag{5.4}
\end{equation}%
\begin{equation*}
\sup\limits_{\lambda \in S\left( \varphi _{1}\right) }R\left\{ \xi ^{\beta
}D^{\beta }\sigma _{2,\lambda }\left( \xi \right) :\xi \in R^{n}/\left\{
0\right\} ,\text{ }\beta \in U_{n}\right\} \leq C.
\end{equation*}%
The estimates $\left( 5.4\right) $ by Theorem 3.1 imply that the operator
functions $\sigma _{1,\lambda }\left( \xi \right) $ and $\sigma _{2,\lambda
}\left( \xi \right) $ are $L_{p,\gamma }\left( R^{n};E\right) $ multipliers.

Let $L_{0}$ denote\ the differential operator in $X$ that\ generated by
problem $\left( 5.1\right) $ for $\lambda =0,$ that is 
\begin{equation*}
D\left( L_{0}\right) =Y,\ L_{0}u=\sum\limits_{\left\vert \alpha \right\vert
=2l}a_{\alpha }D^{\alpha }u+Au.
\end{equation*}%
The estimate $\left( 5.2\right) $ implies that the operator $L_{0}$ has a
bounded inverse from $X$ into $Y.$ We denote by $L$ differential operator in 
$X$ that\ generated by problem $\left( 1.1\right) $, i.e. 
\begin{equation*}
D\left( L\right) =Y,\text{ }Lu=L_{0}u+L_{1}u,\text{ }L_{1}u=\sum\limits_{%
\left\vert \alpha \right\vert \leq 2l}A_{\alpha }\left( x\right) D^{\alpha
}u.
\end{equation*}

\textbf{Theorem 5.2.}\ Suppose all conditions of Theorem 5.1 are hold and 
\begin{equation*}
A_{\alpha }\left( x\right) A^{-\left( 1-\frac{\left\vert \alpha \right\vert 
}{2l}-\mu \right) }\in L_{\infty }\left( R^{n};B\left( E\right) \right) 
\text{ for }0<\mu <1-\frac{\left\vert \alpha \right\vert }{2l}.
\end{equation*}%
Then for all $\ f\in X$ and $\lambda \in S\left( \varphi _{1}\right) $ with
sufficiently large $\left\vert \lambda \right\vert $ equation $\left(
1.1\right) $ has an unique solution $u$ that belongs to space $Y$ and the
uniform coercive estimate holds

\begin{equation}
\sum\limits_{\left\vert \alpha \right\vert \leq 2l}\left\vert \lambda
\right\vert ^{1-\frac{\left\vert \alpha \right\vert }{2l}}\left\Vert
D^{\alpha }u\right\Vert _{X}+\left\Vert Au\right\Vert _{X}\leq C\left\Vert
f\right\Vert _{X}.  \tag{5.5}
\end{equation}

\ \textbf{Proof}. In view of condition on $A_{\alpha }\left( x\right) $ and
by virtue of Theorem 4.1 there is $h>0$ such that \ \ \ 

\begin{equation}
\left\Vert L_{1}u\right\Vert _{X}\leq \sum\limits_{\left\vert \alpha
\right\vert <2l}\left\Vert A_{\alpha }\left( x\right) D^{\alpha
}u\right\Vert _{X}\leq C\sum\limits_{\left\vert \alpha \right\vert
<2l}\left\Vert A^{1-\frac{\left\vert \alpha \right\vert }{2l}-\mu }D^{\alpha
}u\right\Vert _{X}\leq  \tag{5.6}
\end{equation}%
\begin{equation*}
h^{\mu }\left( \ \sum\limits_{\left\vert \alpha \right\vert =2l}\left\Vert
D^{\alpha }u\right\Vert _{X}+\left\Vert \left( A+\lambda \right)
u\right\Vert _{X}\right) +h^{-\left( 1-\mu \right) }\left\Vert u\right\Vert
_{X}
\end{equation*}%
for $u\in Y$. Then from estimates $\left( 5.2\right) $ and $\left(
5.6\right) $ for\ $u\in Y$ we\ have

\begin{equation}
\left\Vert L_{1}u\right\Vert _{X}\leq C\left[ h^{\mu }\left\Vert
(L_{0}+\lambda )u\right\Vert _{X}+h^{-\left( 1-\mu \right) }\left\Vert
u\right\Vert _{X}\right] .  \tag{5.7}
\end{equation}%
Since $\left\Vert u\right\Vert _{X}=\frac{1}{\lambda }\left\Vert \left(
L_{0}+\lambda \right) u-L_{0}u\right\Vert _{X}$ for\ $u\in Y$ we get 
\begin{equation}
\left\Vert u\right\Vert _{X}\leq \frac{1}{\lambda }\left\Vert \left(
L_{0}+\lambda \right) u\right\Vert _{X}+\left\Vert L_{0}u\right\Vert _{X}\leq
\tag{5.8}
\end{equation}%
\begin{equation*}
\frac{1}{\lambda }\left\Vert \left( L_{0}+\lambda \right) u\right\Vert _{X}++%
\frac{M}{\lambda }\left( \sum\limits_{\left\vert \alpha \right\vert
=2l}\left\Vert D^{\alpha }u\right\Vert _{X}+\left\Vert Au\right\Vert
_{X}\right) .
\end{equation*}%
From estimates $\left( 5.7\right) $ and $\left( 5.8\right) $ for $u\in Y$ we
obtain 
\begin{equation}
\left\Vert L_{1}u\right\Vert _{X}\leq Ch^{\mu }\left\Vert \left(
L_{0}+\lambda \right) u\right\Vert _{X}+CM\lambda ^{-1}h^{-\left( 1-\mu
\right) }\left\Vert \left( L_{0}+\lambda \right) u\right\Vert _{X}. 
\tag{5.9}
\end{equation}%
Then choosing $h$ and $\lambda $ such that\ $Ch^{\mu }<1,$ $CMh^{-\left(
1-\mu \right) }<\lambda ,$ from $\left( 5.9\right) $ for sufficiently large $%
\lambda $ we have 
\begin{equation}
\ \ \left\Vert L_{1}\left( L_{0}+\lambda \right) ^{-1}\right\Vert _{B\left(
X\right) }<1.  \tag{5.10}
\end{equation}%
Since we have the relation 
\begin{equation*}
\text{ }\left( L+\lambda \right) ^{-1}=\left( L_{0}+\lambda \right) ^{-1}%
\left[ I+L_{1}\left( L_{0}+\lambda \right) ^{-1}\right] ^{-1}
\end{equation*}%
so by using the estimates $\left( 5.5\right) ,$ $\left( 5.10\right) $ and
the perturbation theory of linear operators we obtain the assertion.

From Theorem 5.2 we obtain the following results:

\textbf{Result 5.1. }Assume the conditions of Theorem 5.2 are satisfied.
Then there exists a constant $C_{1}$ and $C_{2}$ depending only on $p$, $%
\gamma $ such that 
\begin{equation*}
C_{1}\left\Vert u\right\Vert _{Y}\leq \left\Vert \left( L+d\right)
u\right\Vert _{X}\leq C_{2}\left\Vert u\right\Vert _{Y}
\end{equation*}%
for all $u\in Y$ and for sufficiently large $d>0.$

\textbf{Result 5.2. }Assume the conditions of Theorem 5.2 are satisfied.
Then the resolvent operator $\left( L+\lambda \right) ^{-1}$ satisfies the
following coercive sharp estimate holds

\begin{equation*}
\sum\limits_{\left\vert \alpha \right\vert \leq 2l}\left\vert \lambda
\right\vert ^{1-\frac{\left\vert \alpha \right\vert }{2l}}\left\Vert
D^{\alpha }\left( L+\lambda \right) ^{-1}\right\Vert _{B\left( X\right)
}+\left\Vert A\left( L+\lambda \right) ^{-1}\right\Vert _{B\left( X\right)
}\leq C
\end{equation*}%
for $\lambda \in S\left( \varphi _{1}\right) .$

The\textbf{\ }Result 5.2\textbf{\ }implies that operator $L$\textbf{\ }is
positive operator in $X$. Then by virtue of $\left[ \text{27, \S 1.14.5}%
\right] $ the operator $L$ is a generator of an analytic semigroup in $X$
for $\varphi \in \left( \frac{\pi }{2},\pi \right) .$

\begin{center}
\textbf{6. The Cauchy problem for abstract parabolic equation }
\end{center}

Consider now, the Cauchy problem $\left( 1.3\right) .$ In this section we
obta\i n the existence and uniqueness of the maximal regular solution of
problem $\left( 1.3\right) $. First all of we show

\textbf{Theorem 6.1. }Assume the conditions of Theorem 5.1 are satisfied.
Then the operator $L_{0}$ is $R$-positive in $X.$

\textbf{Proof. }Theorem 5.1 implies that the operator $L_{0}$ is positive in 
$X$. We have to prove the $R$-boundedness of the set 
\begin{equation*}
\sigma \left( \lambda \right) =\left\{ \lambda \left( L_{0}+\lambda \right)
^{-1}:\lambda \in S_{\varphi }\right\} .
\end{equation*}%
From Theorem 5.1 we have 
\begin{equation*}
\lambda \left( L_{0}+\lambda \right) ^{-1}f=F^{-1}\Phi \left( \xi ,\lambda
\right) \hat{f}\text{, }
\end{equation*}%
for $f\in X,$ where 
\begin{equation*}
\Phi \left( \xi ,\lambda \right) =\lambda \left( A+L_{0}\left( \xi \right)
+\lambda \right) ^{-1},\text{ }L_{0}\left( \xi \right)
=\sum\limits_{\left\vert \alpha \right\vert =2l}a_{\alpha }\xi ^{\alpha }.
\end{equation*}%
By definition of $R$-boundedness, it is sufficient to show that the operator
function $\Phi \left( \xi ,\lambda \right) $ (depended on variable $\lambda $
and parameters $\xi ,$ $\varepsilon $ ) is uniformly bounded multiplier in $%
X.$ In a similar manner one can easily show that $\Phi \left( \xi ,\lambda
\right) $ is multiplier in $X.$ Then,\ by definition of $R$-boundedness we
have 
\begin{equation*}
\int\limits_{0}^{1}\left\Vert \sum\limits_{j=1}^{m}r_{j}\left( y\right)
\lambda _{j}\left( L_{0}+\lambda _{j}\right) ^{-1}f_{j}\right\Vert
_{X}dy=\int\limits_{0}^{1}\left\Vert \sum\limits_{j=1}^{m}r_{j}\left(
y\right) F^{-1}\Phi \left( \xi ,\lambda _{j}\right) \hat{f}_{j}\right\Vert
_{X}dy=
\end{equation*}

\begin{equation*}
\int\limits_{0}^{1}\left\Vert F^{-1}\sum\limits_{j=1}^{m}r_{j}\left(
y\right) \Phi \left( \xi ,\lambda _{j}\right) \hat{f}_{j}\right\Vert
_{X}dy\leq C\int\limits_{0}^{1}\left\Vert \sum\limits_{j=1}^{m}r_{j}\left(
y\right) f_{j}\right\Vert _{X}dy
\end{equation*}%
for all $\xi _{1},\xi _{2},...,\xi _{m}\in R^{n}$, $\lambda _{1},\lambda
_{2},...,\lambda _{m}\in S_{\varphi },$ $f_{1,}f_{2},...,f_{m}\in X$, $m\in 
\mathbb{N}$, where $\left\{ r_{j}\right\} $ is a sequence of independent
symmetric $\left\{ -1,1\right\} $-valued random variables on $\left[ 0,1%
\right] $. Hence, the set $\sigma \left( \lambda \right) $ is $R$-bounded.

\bigskip For $\mathbf{p=}\left( p,p_{1}\right) $, $R_{+}^{n+1}=R_{+}\times
R^{n},$ $F=L_{\mathbf{p,}\gamma }\left( R_{+}^{n+1};E\right) $ will be
denoted the space of all $E$-valued $\mathbf{p}$-summable weighted functions
with mixed norm, i.e. the space of all measurable functions $f$ defined on $%
R_{+}^{n+1}$ for which 
\begin{equation*}
\left\Vert f\right\Vert _{L_{\mathbf{p,}\gamma }\left( R_{+}^{n+1};E\right)
}=\left( \int\limits_{R_{+}}\left( \dint\limits_{R^{n}}\left\Vert f\left(
x\right) \right\Vert ^{p}\gamma \left( x\right) dx\right) ^{\frac{p_{1}}{p}%
}dt\right) ^{\frac{1}{p_{1}}}<\infty .
\end{equation*}%
Analogously, $F_{0}=W_{\mathbf{p,}\gamma }^{1,2l}\left( R_{+}^{n+1},E\left(
A\right) ,E\right) $ denotes the Sobolev-Lions space with corresponding
mixed norm, i.e.%
\begin{equation*}
F_{0}=\left\{ u\text{: }u\in F\text{,}\right. \text{ }\frac{\partial u}{%
\partial t}\in F\text{, }D^{\alpha }u\in F\text{, }\left\vert \alpha
\right\vert =2l,
\end{equation*}%
\begin{equation*}
\text{ }\left\Vert u\right\Vert _{Y}=\left\Vert \frac{\partial u}{\partial t}%
\right\Vert _{F}+\dsum\limits_{\left\vert \alpha \right\vert =2l}\left\Vert
D^{\alpha }u\right\Vert _{F}+\left\Vert Au\right\Vert _{F}<\infty .
\end{equation*}

The main result of this section is the following:

\textbf{Theorem 6.2.}\ Assume all conditions of Theorem 5.1 hold for $%
\varphi \in \left( \frac{\pi }{2},\pi \right) $ and $p_{1}\in \left(
1,\infty \right) $. Then for $f\in F$ problem $\left( 1.3\right) $ has a
unique solution $u\in F_{0}$ satisfying 
\begin{equation}
\left\Vert \partial _{t}u\right\Vert _{F}+\sum\limits_{\left\vert \alpha
\right\vert =2l}\left\Vert D^{\alpha }u\right\Vert _{F}+\left\Vert
Au\right\Vert _{F}\leq C\left\Vert f\right\Vert _{F}.  \tag{6.1}
\end{equation}

\textbf{Proof.} So, the problem $\left( 1.3\right) $\ can be expressed as 
\begin{equation}
\frac{du}{dt}+L_{0}u\left( t\right) =f\left( t\right) ,\text{ }u\left(
0\right) =0,\text{ }t\in \left( 0,\infty \right) .  \tag{6.2}
\end{equation}

By the Result 5.2 the operator $L_{0}$ is positive in $X$. The Theorem 6.1
implies that $L_{0}$ is $R-$positivity in $X$ for $\varphi \in \left( \frac{%
\pi }{2},\pi \right) .$ Then by virtue of $\left[ \text{29, Th. 4.10}\right] 
$ we obtain that, for $f\in L_{p_{1}}\left( R_{+};X\right) $ the Cauchy
problem $\left( 6.2\right) $ has a unique solution $u\in F_{0}$ satisfying 
\begin{equation}
\left\Vert D_{t}u\right\Vert _{L_{p_{1}}\left( R_{+};X\right) }+\left\Vert
L_{0}u\right\Vert _{L_{p_{1}}\left( R_{+};X\right) }\leq C\left\Vert
f\right\Vert _{L_{p_{1}}\left( R_{+};X\right) }.  \tag{6.3}
\end{equation}

In view of Result 5.1 the operator $L_{0}$ is separable in $X,$ i.e, the
estimate $\left( 6.3\right) $ implies $\left( 6.1\right) $.

\begin{center}
\textbf{\ 7. Degenerate abstract differential equations }
\end{center}

\bigskip\ Let us consider the problem%
\begin{equation}
\ Lu=\sum\limits_{\left\vert \alpha \right\vert =2l}a_{\alpha }D^{\left[
\alpha \right] }u+Au+\sum\limits_{\left\vert \alpha \right\vert
<2l}A_{\alpha }\left( x\right) D^{\left[ \alpha \right] }u+\lambda u=f, 
\tag{7.1}
\end{equation}%
where\ $A$, $A_{\alpha }$\ are linear operators in a Banach space $E$ and $%
\lambda $ is a complex parameter, where%
\begin{equation*}
D_{k}^{\left[ \alpha _{k}\right] }=\left( \gamma _{k}\left( x_{k}\right) 
\frac{\partial }{\partial x_{k}}\right) ^{\alpha _{k}}\text{, }D^{\left[
\alpha \right] }=D_{1}^{\left[ \alpha _{1}\right] }D_{2}^{\left[ \alpha _{2}%
\right] }...D_{n}^{\left[ \alpha _{n}\right] },
\end{equation*}%
here $\gamma _{k}\left( x\right) $ are positive measurable functions on $%
R^{n}.$

Let

\begin{equation*}
W_{p\mathbf{,}\gamma }^{\left[ l\right] }\left( \Omega ,E_{0},E\right)
=\left\{ u\in L_{p}\left( \Omega ;E_{0}\right) ,\text{ }D_{k}^{\left[ l_{k}%
\right] }u\in L_{p}\left( \Omega ;E\right) \right\} ,
\end{equation*}%
\begin{equation*}
\left\Vert u\right\Vert _{W_{p\mathbf{,\gamma }}^{\left[ l\right] }\left(
\Omega ;E_{0,}E\right) }=\left\Vert u\right\Vert _{L_{p}\left( \Omega
;E_{0}\right) }+\sum\limits_{k=1}^{n}\left\Vert D_{k}^{\left[ l_{k}\right]
}u\right\Vert _{L_{p}\left( \Omega ;E\right) }<\infty .
\end{equation*}

Here, 
\begin{equation*}
X=L_{p}\left( R^{n};E\right) ,\text{ }Y=W_{p\mathbf{,}\gamma }^{\left[ 2l%
\right] }\left( R^{n};E\left( A\right) ,E\right) .
\end{equation*}

Let 
\begin{equation}
\int\limits_{0}^{x_{k}}\gamma _{k}^{-1}\left( y\right) dy<\infty \text{, }%
k=1,2,...,n.  \tag{7.2}
\end{equation}

\textbf{Remark 7.1. }\bigskip Under the substitution 
\begin{equation}
\tau _{k}=\int\limits_{0}^{x_{k}}\gamma _{k}^{-1}\left( y\right) dy 
\tag{7.3}
\end{equation}%
the spaces $X$ and $Y$\ are mapped isomorphically onto the weighted spaces $%
L_{p,\tilde{\gamma}}\left( R^{n};E\right) $, $W_{p,\tilde{\gamma}%
}^{2l}\left( R^{n};E\left( A\right) ,E\right) ,$ where%
\begin{equation*}
\tilde{\gamma}=\tilde{\gamma}\left( \tau \right)
=\prod\limits_{k=1}^{n}\gamma _{k}\left( x_{k}\left( \tau _{k}\right)
\right) \text{, }\tau =\left( \tau _{1},\tau _{2},...,\tau _{n}\right) .
\end{equation*}

Moreover, under the transformation $\left( 7.3\right) $ the problem $\left(
7.1\right) $ is mapped to the undegenerate problem $\left( 1.1\right) $
considered in the weighted space $L_{p,\tilde{\gamma}}\left( R^{n};E\right) $%
.

\textbf{Condition 7.1. }Assume $\left( 7.1\right) $ holds and $\gamma
_{k}\left( x_{k}\left( \tau _{k}\right) \right) \in A_{p}$ for $k=1,2,...,n$
and $p\in \left( 1,\infty \right) .$

From Theorem 5.2 and Remark 7.1 we obtain the following results:

\textbf{Result 7.1. }Assume the conditions of Theorem 5.2 are satisfied.
Then for all $\ f\in X$ and $\lambda \in S\left( \varphi _{1}\right) $ with
sufficiently large $\left\vert \lambda \right\vert $ equation $\left(
1.1\right) $ has an unique solution $u$ that belongs to $Y$ and the uniform
coercive estimate holds

\begin{equation*}
\sum\limits_{\left\vert \alpha \right\vert \leq 2l}\left\vert \lambda
\right\vert ^{1-\frac{\left\vert \alpha \right\vert }{2l}}\left\Vert D^{%
\left[ \alpha \right] }u\right\Vert _{X}+\left\Vert Au\right\Vert _{X}\leq
C\left\Vert f\right\Vert _{X}.
\end{equation*}

Let $G$ denote the operator in $X$ generated by the problem $\left(
7.1\right) .$

\textbf{Result 7.2. }Assume the conditions of Theorem 5.2 and the Condition
7.1 are satisfied. Then the resolvent operator $\left( L+\lambda \right)
^{-1}$ satisfies the following sharp estimate

\begin{equation*}
\sum\limits_{\left\vert \alpha \right\vert \leq 2l}\left\vert \lambda
\right\vert ^{1-\frac{\left\vert \alpha \right\vert }{2l}}\left\Vert D^{%
\left[ \alpha \right] }\left( G+\lambda \right) ^{-1}\right\Vert _{B\left(
X\right) }+\left\Vert A\left( G+\lambda \right) ^{-1}\right\Vert _{B\left(
X\right) }\leq C
\end{equation*}%
for $\lambda \in S\left( \varphi _{1}\right) .$

The\textbf{\ }Result 5.2\textbf{\ }implies that operator $G$\textbf{\ }is
positive operator in $X$. Then by virtue of $\left[ \text{27, \S 1.14.5}%
\right] $ the operator $G$ is a generator of an analytic semigroup in $X$
for $\varphi \in \left( \frac{\pi }{2},\pi \right) .$

Consider the Cauchy problem for degenerate parabolic equation%
\begin{equation}
\ \partial _{t}u+\sum\limits_{\left\vert \alpha \right\vert =2l}a_{\alpha
}D^{\left[ \alpha \right] }u+Au=f\left( t,x\right) ,\text{ }t\in \left(
0,\infty \right) \text{, }x\in R^{n},  \tag{7.4}
\end{equation}%
\begin{equation*}
u\left( 0,x\right) =0,\text{ }x\in R^{n},
\end{equation*}%
where $a_{\alpha }$ are complex numbers and $A$ is a linear operator in a
Banach space $E.$

For $\mathbf{p=}\left( p,p_{1}\right) $, let $\Phi =L_{\mathbf{p}}\left(
R_{+}^{n+1};E\right) $ denotes $L_{\mathbf{p,}\gamma }\left(
R_{+}^{n+1};E\right) $ for $\gamma \left( x\right) \equiv 1.$ Analogously, $%
\Phi _{0}=W_{\mathbf{p,}\gamma }^{1,\left[ 2l\right] }\left(
R_{+}^{n+1},E\left( A\right) ,E\right) $ denotes the Sobolev-Lions space
with corresponding mixed norm, i.e.%
\begin{equation*}
\Phi _{0}=\left\{ u\text{: }u\in \Phi \text{,}\right. \text{ }\frac{\partial
u}{\partial t}\in \Phi \text{, }D^{\left[ \alpha \right] }u\in \Phi \text{, }%
\left\vert \alpha \right\vert =2l,
\end{equation*}%
\begin{equation*}
\text{ }\left\Vert u\right\Vert _{\Phi _{0}}=\left\Vert \frac{\partial u}{%
\partial t}\right\Vert _{\Phi }+\dsum\limits_{\left\vert \alpha \right\vert
=2l}\left\Vert D^{\left[ \alpha \right] }u\right\Vert _{\Phi }+\left\Vert
Au\right\Vert _{\Phi }<\infty .
\end{equation*}

From Theorem 6.2 and Remark 7.1 we obtain the following results:

\bigskip \textbf{Result 7.3. \ }Assume all conditions of Theorem 5.1 and the
Condition 7.1 are satisfied for $\varphi \in \left( \frac{\pi }{2},\pi
\right) $ and $p_{1}\in \left( 1,\infty \right) $. Then for all $f\in \Phi $
problem $\left( 7.4\right) $ has a unique solution $u\in \Phi _{0}$
satisfying 
\begin{equation*}
\left\Vert \frac{\partial u}{\partial t}\right\Vert _{\Phi
}+\sum\limits_{\left\vert \alpha \right\vert =2l}\left\Vert D^{\left[ \alpha %
\right] }u\right\Vert _{\Phi }+\left\Vert Au\right\Vert _{\Phi }\leq
C\left\Vert f\right\Vert _{\Phi }.
\end{equation*}

\begin{center}
\textbf{8.} \textbf{Maximal regularity properties of infinite many system of
parabolic equations}
\end{center}

Consider the Cauchy problem for infinite many system of parabolic equations

\begin{equation}
\partial _{t}u_{i}\left( t,x\right) \sum\limits_{\left\vert \alpha
\right\vert =2l}a_{\alpha }D^{\alpha }u_{i}\left( t,x\right)
+\sum\limits_{j=1}^{\infty }a_{ij}u_{j}\left( t,x\right) =f_{i}\left(
t,x\right) \text{, }x\in R^{n}\text{, }t\in \left( 0,\infty \right) , 
\tag{8.1}
\end{equation}

\begin{equation}
u\left( 0,x\right) =0,\text{ for a.e. }x\in R^{n}\text{, }i=1,2,...,N,\text{ 
}N\in \mathbb{N},  \tag{8.2}
\end{equation}%
where $a_{\alpha }$ and $a_{ij}$ are complex numbers.

\textbf{Condition 8.1. }Let 
\begin{equation*}
a_{ij}=a_{ji}\text{, }\sum\limits_{i,j=1}^{N}a_{ij}\xi _{i}\xi _{j}\geq
C_{0}\left\vert \xi \right\vert ^{2},\text{ for }\xi \neq 0.
\end{equation*}

\bigskip Let%
\begin{equation*}
\text{ }u=\left\{ u_{j}\right\} ,\text{ }Au=\left\{
\dsum\limits_{j=1}^{N}a_{ij}u_{j}\right\} ,\text{ }i,\text{ }j=1,2,...N,
\end{equation*}
\begin{equation*}
\text{ }l_{q}\left( A\right) =\left\{ u\in l_{q},\left\Vert u\right\Vert
_{l_{q}\left( A\right) }=\left\Vert Au\right\Vert _{l_{q}}=\right.
\end{equation*}

\begin{equation*}
\left. \left( \sum\limits_{i=1}^{N}\left\vert
\sum\limits_{j=1}^{N}a_{ij}u_{j}\right\vert ^{q}\right) ^{\frac{1}{q}%
}<\infty \right\} ,\text{ }q\in \left( 1,\infty \right) .
\end{equation*}

Here, 
\begin{equation*}
X_{\mathbf{p},q}=L_{\mathbf{p},\gamma }\left( R^{n};l_{q}\right) ,\text{ }Y_{%
\mathbf{p},q}=W_{\mathbf{p},\gamma }^{1,2l}\left( R_{+}^{n+1},l_{q}\left(
A\right) ,l_{q}\right)
\end{equation*}

\textbf{Theorem 8.1.} Assume the Conditions 5.1 and 8.1 are satisfied. Then
for all $f\left( x\right) =\left\{ f_{i}\left( x\right) \right\}
_{1}^{\infty }\in X_{\mathbf{p},q}$ problem $\left( 8.1\right) -\left(
8.2\right) $ has a unique solution $u=\left\{ u_{i}\left( x\right) \right\}
_{1}^{\infty }$ that belongs to space $Y_{\mathbf{p},q}$ and the coercive
sharp estimate holds%
\begin{equation}
\left\Vert \frac{\partial u}{\partial t}\right\Vert _{X_{\mathbf{p}%
,q}}+\sum\limits_{\left\vert \alpha \right\vert =2l}\left\Vert D^{\alpha
}u\right\Vert _{X_{\mathbf{p},q}}+\left\Vert Au\right\Vert _{X_{\mathbf{p}%
,q}}\leq C\left\Vert f\right\Vert _{X_{\mathbf{p},q}}.  \tag{8.3}
\end{equation}%
\ \textbf{Proof. }Let $E=l_{q},$ $A$ be a matrix such that $A=\left[ a_{ij}%
\right] ,$ $i$, $j=1,2,...N.$ It is easy to see that 
\begin{equation*}
B\left( \lambda \right) =\lambda \left( A+\lambda \right) ^{-1}=\frac{%
\lambda }{D\left( \lambda \right) }\left[ A_{ji}\left( \lambda \right) %
\right] \text{, }i\text{, }j=1,2,...N,
\end{equation*}%
where $D\left( \lambda \right) =\det \left( A-\lambda I\right) $, $%
A_{ji}\left( \lambda \right) $ are entries of the corresponding adjoint
matrix of $A-\lambda I.$ Since the matrix $A$ is symmetric and positive
definite, it generates a positive operator in $l_{q}$ for $q\in \left(
1,\infty \right) .$ For all $u_{1,}u_{2},...,u_{\mu }\in l_{q}$, $\lambda
_{1},\lambda _{2},...,\lambda _{\mu }\in \mathbb{C}$ and independent
symmetric $\left\{ -1,1\right\} $-valued random variables $r_{k}\left(
y\right) $, $k=1,2,...,\mu ,$ $\mu \in \mathbb{N}$\ we have 
\begin{equation*}
\int\limits_{\Omega }\left\Vert \sum\limits_{k=1}^{\mu }r_{k}\left( y\right)
B\left( \lambda _{k}\right) u_{k}\right\Vert _{l_{q}}^{q}dy\leq
\end{equation*}%
\begin{equation*}
C\left\{ \int\limits_{\Omega }\sum\limits_{j=1}^{N}\left\vert
\sum\limits_{k=1}^{\mu }\sum\limits_{j=1}^{N}\frac{\lambda _{k}}{D\left(
\lambda _{k}\right) }A_{ji}\left( \lambda _{k}\right) r_{k}\left( y\right)
u_{ki}\right\vert ^{q}dy\right. \leq
\end{equation*}%
\begin{equation}
\sup\limits_{k,i}\sum\limits_{j=1}^{N}\left\vert \frac{\lambda _{k}}{D\left(
\lambda _{k}\right) }A_{ji}\left( \lambda _{k}\right) \right\vert
^{q}\int\limits_{\Omega }\left\vert \sum\limits_{k=1}^{\mu }r_{k}\left(
y\right) u_{kj}\right\vert ^{q}dy.  \tag{8.4}
\end{equation}%
Since $A$ is symmetric and positive definite, we have%
\begin{equation}
\sup\limits_{k,i}\sum\limits_{j=1}^{N}\left\vert \frac{\lambda _{k}}{D\left(
\lambda _{k}\right) }A_{ji}\left( \lambda _{k}\right) \right\vert ^{q}\leq C.
\tag{8.5}
\end{equation}%
From $\left( 8.4\right) $ and $\left( 8.5\right) $ we get 
\begin{equation*}
\int\limits_{\Omega }\left\Vert \sum\limits_{k=1}^{\mu }r_{k}\left( y\right)
B\left( \lambda _{k}\right) u_{k}\right\Vert _{l_{q}}^{q}dy\leq
C\int\limits_{\Omega }\left\Vert \sum\limits_{k=1}^{\mu }r_{k}\left(
y\right) u_{k}\right\Vert _{l_{q}}^{q}dy.
\end{equation*}%
i.e., the operator $A$ is $R$-positive in $l_{q}.$ Hence, by Theorem 6.2 we
obtain the assertion.

\ \textbf{Remark 8.1}. There are a lot of $R-$positive operators in
different concrete Banach spaces. Therefore, putting concrete Banach spaces
instead of $E,$ and concrete differential, pseudo differential operators, or
finite, infinite matrices instead of $A,$ by virtue of Theorems 5.2 \ and
6.2 we can obtained the different class of maximal regular partial
differential equations or system of equations.

\begin{center}
\bigskip \textbf{Acknowledgements}
\end{center}

The author would like to express a gratitude to Dr. Neil. Course for his
useful advice in English in preparing of this paper

\bigskip

\begin{center}
\ \textbf{References}
\end{center}

\bigskip\ \ \ \ \ \ \ \ \ \ \ \ \ \ \ \ \ \ \ \ \ \ \ \ \ \ \ \ \ \ \ \ \ \
\ \ \ \ \ \ \ \ \ \ \ \ \ \ \ \ \ \ \ \ \ \ \ \ \ \ \ \ \ \ \ \ \ \ \ \ \ \
\ \ \ \ \ \ \ \ \ \ \ \ \ \ \ \ \ \ \ 

\begin{enumerate}
\item Amann, H., Linear and quasi-linear equations,1, Birkhauser, Basel 1995.

\item Agarwal, R., O' Regan, D., Shakhmurov V. B., Separable anisotropic
differential operators in weighted abstract spaces and applications, J.
Math. Anal. Appl. 338(2008), 970-983.

\item Besov, O. V., Ilin, V. P., Nikolskii, S. M., Integral representations
of functions and embedding theorems, Nauka, Moscow, 1975.

\item Bourgain, J., Vector--Valued Singular Integrals and the H1--BMO
Duality. In: Probability theory and harmonic analysis, pp. 1 -- 19, Pure
Appl. Math. 98, Marcel Dekker, 1986.

\item Cl\'{e}ment Ph., De Pagter B., Sukochev F. A., Witvliet H., Schauder
decomposition and multiplier theorems, Studia Math. 138 (2000), 135-163.

\item Chill, R., Fiorenza, A., Singular integral operators with
operator-valued kernels, and extrapolation of maximal regularity into
rearrangement invariant Banach function spaces, J. Evol. Equ. 14 (2014),
4-5, 795--828.

\item Cl\'{e}ment, P., Pr\"{u}ss, J., An operator--valued transference
principle and maximal regularity on vector--valued Lp--spaces, Evolutian
equations and their applications in physical and life sciences, p. 67 -- 87,
Bad Herrenalb, 1998.

\item Denk, R., Hieber M., Pr\"{u}ss J., $R-$boundedness, Fourier
multipliers and problems of elliptic and parabolic type, Mem. Amer. Math.
Soc. 166 (2003), n.788.

\item Haller, R., Heck H., Noll A., Mikhlin's theorem for operator-valued
Fourier multipliers in $n$ variables, Math. Nachr. 244, (2002), 110-130.

\item Girardi, M., Lutz, W., Operator-valued Fourier multiplier theorems on $%
L_{p}$($X$) and geometry of Banach spaces, J. Funct. Anal., 204(2),
320--354, 2003.

\item Hyt\"{o}nen, T. P., Anisotropic Fourier multipliers and singular
integrals for vector-valued functions, Ann. Mat. Pura Appl. (4) 186
(2007)(3), 455--468.

\item H\"{a}nninen, T. S., Hyt\"{o}nen, T. P., The A2 theorem and the local
oscillation decomposition for Banach space valued functions, J. Operator
Theory, 72 (2014), (1), 193--218.

\item Kurtz, D. G., Whedeen, R. L., Results on weighted norm inequalities
for multipliers, Trans. Amer. Math. Soc. 255(1979), 343--362.

\item Lions, J. L and Peetre J., Sur une classe d'espaces d'interpolation,
Inst. Hautes Etudes Sci. Publ. Math., 19(1964), 5-68.

\item Lizorkin, P. I, Shakhmurov, V. B., Embedding theorem for vector-valued
functions. 1., Izvestiya Vushkh Uchebnykh Zavedenie Matematika (1)1989,
70-79; 2. (8) 1989, 69-78.

\item Meyries, M., Veraar, M., Pointwise multiplication on vector-valued
function spaces with power weights. J. Fourier Anal. Appl. 21 (2015)(1),
95--136.

\item Meyries M. and Veraar, M. C., Sharp embedding results for spaces of
smooth functions with power weights, Studia Math., 208(3):257--293, 2012.

\item McConnell Terry R., On Fourier Multiplier Transformations of
Banach-Valued Functions, Trans. Amer. Mat. Soc. 285, (2) (1984), 739-757.

\item Pisier, G., Some Results on Banach Spaces without Local Unconditional
Structure, Compositio Math. 37 (1978), 3 -- 19.

\item Ragusa, M. A., Embeddings for Lorentz-Morrey spaces, J. Optim. Theory
Appl., 154(2)(2012), 491-499.

\item Shakhmurov V. B, Abstract capacity of regions and compact embedding
with applications, Acta. Math. Scia., (31)1, 2011, 49-67.

\item Shakhmurov, V. B., Imbedding theorems and their applications to
degenerate equations, Differential equations, 24 (4), (1988), 475-482.

\item Shakhmurov, V. B., Embedding operators and maximal regular
differential-operator equations in Banach-valued function spaces, J.
Inequal. Appl., 4(2005), 329-345.

\item Shakhmurov, V. B., Coercive boundary value problems for regular
degenerate differential-operator equations, J. Math. Anal. Appl., 292 ( 2),
(2004), 605-620.

\item Shakhmurov, V. B., Embedding and maximal regular differential
operators in Sobolev-Lions spaces, Acta Mathematica Sinica, 22(5) 2006,
1493-1508.

\item Schmeisser H., Vector-valued Sobolev and Besov spaces, Seminar
analysis of the Karl-Weierstra -Institute of Mathematics 1985/86 (Berlin,
1985/86), Teubner, Leipzig, 1987, 4-44.

\item Triebel, H., Interpolation theory. Function spaces. Differential
operators, North-Holland, Amsterdam, 1978.

\item Triebel, H., Spaces of distributions with weights. Multipliers in $%
L_{p}$-spaces with weights, Math. Nachr. 78, (1977), 339-356.

\item Weis, L., Operator-valued Fourier multiplier theorems and maximal $%
L_{p}$ regularity, Math. Ann. 319, (2001), 735-75.

\item Yakubov, S and Yakubov Ya., Differential-operator equations. Ordinary
and Partial \ Differential equations, Chapman and Hall /CRC, Boca Raton,
2000.

\item Zimmerman, F., On vector-valued Fourier multiplier theorems, Studia
Math. 93 (3)(1989), 201-222.
\end{enumerate}

\bigskip

\begin{center}
\bigskip

\bigskip
\end{center}

\ \ \ \ \ \ \ \ \ \ \ \ \ \ \ \ \ \ \ \ \ \ \ \ 

\ 

\begin{center}
\bigskip
\end{center}

\end{document}